\documentclass{IEEEtran}
\usepackage{graphicx}          

\usepackage[T1]{fontenc} 
\usepackage{algorithm, algorithmicx, algpseudocode}
\usepackage{amsmath, amsfonts, amssymb, mathrsfs}
\usepackage{fancyhdr}
\usepackage{enumerate}
\usepackage{cite}
\usepackage{epstopdf}
\usepackage{xcolor}

\newtheorem{theorem}{Theorem}
\newtheorem{lemma}{Lemma}
\newtheorem{definition}{Definition}
\newtheorem{assumption}{Assumption}
\newtheorem{remark}{Remark}

\newtheorem{coro}{Corollary}
\newtheorem{example}{Example}
\newenvironment{pf}{\begin{IEEEproof}}{\end{IEEEproof}}

\def\BibTeX{{\rm B\kern-.05em{\sc i\kern-.025em b}\kern-.08em
    T\kern-.1667em\lower.7ex\hbox{E}\kern-.125emX}}
    \def \qed {\hfill \vrule height6pt width 6pt depth 0pt}
\begin{document}
\title{Reliably Learn to Trim Multiparametric Quadratic Programs via Constraint Removal}
\author{Zhinan Hou, Keyou You
\thanks{This work was supported by National Science and Technology Major Project of China (2022ZD0116700) and National Natural Science Foundation of China (62033006, 62325305) (Corresponding author: Keyou You).}
\thanks{Zhinan Hou and Keyou You are with the Department of Automation and Beijing National Research Center for Information Science and Technology, Tsinghua University, Beijing 100084, China. (e-mail: hzn22@mails.tsinghua.edu.cn, youky@tsinghua.edu.cn)}}
\maketitle
\begin{abstract}

In a wide range of applications, we are required to rapidly solve a sequence of  multiparametric quadratic program (mp-QP) instances  on resource-limited hardware. This has been an active topic for decades in control and optimization communities. Observe that the main computational cost of existing algorithms lies in addressing many linear inequality constraints, though their majority are redundant and removing them will not change the optimal solution. This work learns from previously solved mp‑QP instances to develop a novel constraint removal method, which reliably trims unsolved instances and significantly reduces solution cost. Then, we extend to trim the mp-QP instances of model predictive control (MPC) whose parameter vectors are sampled from linear systems. Importantly, both online and offline solved instances can be utilized  to adaptively trim instances in the closed-loop system. We show that the number of linear inequalities in the trimmed instance of MPC decreases to zero in a finite timestep, which can also be reduced by increasing offline computation. Finally, simulations are performed to demonstrate the efficiency of our trimming method in removing redundant constraints.
\end{abstract}


\section{Introduction}

This work focuses on the multiparametric quadratic program (mp-QP) with a convex quadratic cost function and many linear inequality constraints where the linear term of the cost function and right-hand side of the constraints depend on a vector of parameters that may change from one step to another.  It has emerged in a huge variety of applications in control \cite{bemporad2002explicit,darby2007parametric}, operations research \cite{jayasekara2023solving}, machine learning \cite{TSO2020106902} and many others, which require rapidly solving mp-QP instances over resource-limited hardware. 
%
%

By exploiting the unique structure of mp-QP, quite a few innovative methods have been proposed in the literature. They can be roughly categorized as active-set \cite{bartlett2006qpschur, arnstrom2022dual}, interior-point \cite{FRISON20206563}, first-order \cite{hermans2019qpalm}, and nonnegative least squares \cite{bemporad2015quadratic} methods. As a first-order method, the alternating direction method of multipliers (ADMM) \cite{boyd2011distributed} is core to the prominently featured solvers -- OSQP \cite{stellato2020osqp} and TinyMPC \cite{nguyen2024tinympc}, which are developed for application of  the embedded model predictive control (MPC).  All these methods consume a large portion of computational resources in handling inequality constraints. In particular, the worst-case complexity of an active-set method grows exponentially with the number of  constraints \cite{wright2006numerical}. However, the majority of inequality constraints in mp-QP turn out to be redundant in the sense that removing them will not change the optimal solution. On the contrary, if the mp-QP has no inequality constraint, the optimal solution simply takes a linear form of the parameter vector and moreover, the gain remains constant over the space of parameter vector. In fact, the less the number of  inequality constraints in mp-QP, the faster these solution methods converge. This motivates us to trim mp-QP via removing redundant inequalities, after which we only need to solve ``slim'' mp-QP instances, hoping to accelerate the above-mentioned solution methods.

Then, a natural question is how to cheaply and rapidly trim mp-QP, {\em as slim as possible}? Noting that the minimizer of mp-QP is Lipschitz continuous in parameter vectors \cite[Proposition 7.13]{rawlings2017model}, a  solved instance can be potentially exploited to accelerate the unsolved one, especially when their parameter vectors are close. This observation forms the basis of warm start in the vast body of literature \cite{stellato2020osqp}.  However, it usually requires sophisticated design and many solution methods, e.g., interior-point methods, are difficult to warm start. Thus, it is mostly used in practice and usually lacks theoretical guarantee. 

Quite differently, we utilize the results of a solved mp-QP instance to remove redundant inequalities of the unsolved instance. More importantly, we propose a novel method to check the redundancy, which at most  involves twice of scalar number comparisons per inequality constraint. Clearly, such a simple method can be efficiently implemented  on resource-limited hardware, which is essential to the design of microcontrollers. In addition, we quantify the importance of the solved instance in removing redundant constraints. Specifically, we evaluate the closeness of the two instances via classifying the distance of their parameter vectors against a sequence of nonnegative numbers, which can be solved in prior and offline by mixed integer linear programs (MILPs). Under the so-called linear independence constraint qualification (LICQ) \cite{borrelli2017predictive}, we prove that if their parameter distance is less than the above-mentioned positive number, then the number of linear inequalities remaining in the trimmed mp-QP instance can be explicitly upper bounded. It is worth mentioning that this number can be very conservative and the empirical results are confirmed much better than the theory dictated. Moreover, we use our trimming method to learn {\em multiple} solved instances and establish a ``the more the less'' result, i.e., the more the solved instances learnt, the less the inequality constraints in the trimmed instance.  

Then, we extend to {\em adaptively} trim mp-QP instances of MPC \cite[Chapter 2]{rawlings2017model} in the online closed-loop system where the parameter vectors are encountered system state vectors. Since MPC is usually used to synthesize control law for constrained systems, the parameter vectors may converge to the origin under standard design choices. We leverage this result to achieve that after a quantified finite number of timestep, the trimmed mp-QP instances have no constraint. Thereafter, the resultant MPC law simply takes a linear form with a fixed control gain. 

Noticeably, our method also allows us to jointly learn {\em offline} solved instances to trim new instances of mp-QP in the closed-loop system of MPC. We provide a computational tradeoff between the explicit MPC \cite{alessio2009survey} and the conventional implicit MPC that yields a control action for each encountered state vector. To elaborate it, the control law of the explicit MPC is offline obtained from solving mp-QP and implemented online via a table-lookup method, while the control action of the implicit MPC is online iteratively computed from mp-QP instances.   Neither of them can benefit from the results of offline and online solved mp-QP instances. In comparison, our trimmed MPC can easily take advantages of both solved instances, and achieve a freedom between offline and online computation. 

It is worth stressing that removing redundant constraints is not a new technique to accelerate optimization problems, see  \cite{romao2022exact, chinneck1997finding, berbee1987hit, choi2007removing}.  For instance, the redundant constraint defined in \cite{boyd2004convex} refers to constraints whose removal does not change the feasible set. Such constraints are typically removed offline before solving. However, offline removal is often insufficient to accelerate real-time applications like MPC, as it can only remove constraints redundant across all feasible states. Thus, some online methods have recently been developed to accelerate mp-QP instances in MPC via system-theoretic approaches \cite{jost2013accelerating, jost2015online, nouwens2023constraint}, where they heavily exploit system properties, such as  region of activity \cite{jost2013accelerating},  reachability \cite{nouwens2023constraint}, and contractility of the cost function \cite{jost2015online}, to achieve constraint removal. Thus, their methods are limited to instances of MPC for linear dynamical systems. Differently, we exploit previously solved instances via a learning-theoretic approach, leading to some salient features: (a) our method can be easily applied  to remove redundant constraints of instances whose parameter vectors are {\em arbitrarily} generated. (b) The number of removed redundant constraints can be quantified. In particular,  when applied to trim mp-QP instances in the linear closed-loop system of MPC, the number of linear inequalities is proved to decrease to zero, which is not the case for the system-theoretic approach in \cite{jost2013accelerating, jost2015online, nouwens2023constraint}. In fact, they are unable to achieve any guarantee and can only illustrate this phenomenon‌  via simulation. Moreover, simulation results demonstrate that our method is much more efficient in removing redundant inequality constraints. 

The rest of paper is organized as follows. In Section \ref{sec:background}, we introduce the convex mp-QP and our objective. In Section \ref{sec:online}, we reliably learn from the solved mp-QP instances to trim unsolved instances via removing redundant inequality constraints. In Section \ref{sec:offline}, we provide a method to compute the global Lipschitz constant of mp-QP and evaluate the performance of our trimming method. In Section \ref{sec:bridge}, we extend to trim mp-QP of linear MPC. Numerical results are included in Section \ref{sec:sim} and we draw some concluding remarks in Section \ref{sec:conclusion}. 

\textbf{Notation.}  Let $\mathbb{N}_{[a,b]}$ be the set of integers lying in the interval $[a, b]$, i.e., $\mathbb{N}_{[a,b]}= \{ n \in \mathbb{N} | a \le n \le b\}$. Given a matrix $A \in \mathcal{R}^{m \times n}$, let $A_j$ denote its $j$-th row and so is the vector. We say that $M$ is a set-valued mapping from $\mathcal{X}$ to $\mathcal{Y}$, denoted by $M: \mathcal{X} \rightrightarrows \mathcal{Y}$, if for every $x \in \mathcal{X}, \ M(x) \subseteq \mathcal{Y}$.  For a vector (matrix) $X$, we denote its $\ell_2$, $\ell_1$ and $\ell_\infty$ norms by $\Vert X \Vert$, $\Vert X \Vert_1$ and $\Vert X \Vert_\infty$, respectively. Let $\mathcal{B}(q,r)$ be a ball centered at $q \in \mathcal{R}^{n}$ with a radius $r>0$, i.e., $\mathcal{B}(q,r) = \{v \in \mathcal{R}^n |~\Vert v - q\Vert \le r  \}$. If every element of $X$ is positive, we simply write $X>0$. For a real number $x$, the ceiling function $\lceil x \rceil$ returns the smallest integer greater than or equal to $x$.  For two sets $\mathcal{A}$ and $\mathcal{B}$, we define $\mathcal{A}-\mathcal{B} = \{ x |x \in \mathcal{A} , x\notin \mathcal{B}\}$.  
\section{Problem formulation and our objective} \label{sec:background}
In this section, we first introduce the convex multiparametric quadratic program (mp-QP). Then, we outline the objective of learning from solved mp-QP instances to trim (unsolved) instances via removing their redundant inequality constraints, after which the trimmed version can be much cheaper to solve. 
In this work, we focus on the following mp-QP
\begin{subequations}\label{2_5}
	\begin{align}
		\text{QP}(x): \ &  \mathop{\text{min.}}\limits_z \ V(z):=\frac{1}{2}z^THz + x^TFz, \label{2_5_1}\\ 
		\text{s.t.} \ & z \in \mathcal{Z}(x) := \{z \in \mathcal{R}^{n_z}| \ Gz \le  Sx+ w \}, \label{2_5_2}
	\end{align} 
\end{subequations}
where $x \in \mathcal{R}^{n_x}$ is the parameter vector and $z \in \mathcal{R}^{n_z}$ is the decision vector. The objective function  is defined by a positive definite matrix $H \in \mathcal{R}^{n_z \times n_z}$ and $F \in \mathcal{R}^{n_x \times n_z}$ in a quadratic form, and the  feasible set $\mathcal{Z}(x)$  is specified by constant $G \in \mathcal{R}^{n_c \times n_z}, w \in \mathcal{R}^{n_c}$ and $S \in \mathcal{R}^{n_c \times n_x}$ in a linear form. Throughout this work, $\mathcal{Z}(x)$ is a non-empty set. Note that some constraints in \eqref{2_5_2} can be removed without changing the optimal solution.

\begin{definition}[Redundant constraint] We say that a constraint is redundant if
removing it does not change the optimal solution of the mp-QP of \eqref{2_5}.
\end{definition}

In Chapter 4.1.1 \cite{boyd2004convex}, a constraint is said to be redundant if deleting it does not change the feasible set.  This is a stronger definition and usually its number is very limited for a well-designed problem with mp-QP. 

Removing these redundant constraints leads to a ``slim'' mp-QP, which can be solved more cheaply and rapidly. Such an idea has been rarely exploited and motivates the study of trimming mp-QP of this work.

Since the difference between two mp-QP instances in the form of \eqref{2_5}  lies in their parameter vectors, the result of one solved instance can be potentially used to accelerate the other one. A naive and common idea is to use the solved one to {\em warm-start} the iteration process of the other, hoping to reduce the iteration number of optimization algorithms. Though empirically appreciated, it generically lacks theoretical guarantee and from the worst-case point of view, it does not help at all. 

This work takes a quite different perspective to utilize the results of solved instances. Notice that if $\mathcal{Z}(x)$ of the mp-QP contains a large number of linear inequalities, which is indeed the case in many applications, the majority of them are generically redundant. Removing them clearly results in a ``slim'' mp-QP problem which, more importantly, could be much cheaper to solve. For example, the optimal solution of the mp-QP in \eqref{2_5} can be explicitly expressed as $z^*(x)=-H^{-1}F^Tx$ if there is no linear inequality constraint. Then, a follow-up problem is how to cheaply and reliably remove redundant linear inequalities in $\mathcal{Z}(x)$ {\em as many as possible}? 

In this work, we propose a novel method to reliably learn from the results of solved mp-QP instances in the form of \eqref{2_5} to trim new instances via constraint removal. It is worth stressing that though the trimmed mp-QP of this work can be made very slim, it has the same optimal solution as the original one.  Moreover, we extend to solve the linear quadratic MPC \cite[Chapter 7]{rawlings2017model} in the form of mp-QP. 

\section{Reliably learn to trim mp-QP via constraint removal
} \label{sec:online}
In this section, we learn from the results of solved mp-QP instances to trim unsolved instances via constraint removal, without changing the optimal solution. First, we introduce a general rule on the optimality of the trimmed mp-QP. Then, we use the results of solved instances to explicitly design an efficient inequality removal rule by leveraging the Lipschitz continuity of the minimizer of the mp-QP. 

\subsection{The trimmed mp-QP}
To formalize our idea, we consider the following trimmed version of \eqref{2_5} 
\begin{subequations}\label{2_6}
	\begin{align}
		\mathop{\text{min.}}\limits_z \ & \frac{1}{2}z^THz + x^TFz, \label{2_6_1}\\ 
		\text{s.t.} \ & z \in \mathcal{Z}(x, \mathbb{I}(x))=\bigcap\nolimits_{j\in \mathbb{I}(x) } \mathcal{Z}_j(x),  \label{2_6_2}
	\end{align} 
\end{subequations}
where $\mathcal{Z}_j(x) := \{z\in\mathcal{R}^{n_z} \ | \ G_jz \le S_jx+ w_j\}$ and $\mathbb{I}: \mathcal{R}^{n_x} \rightrightarrows \mathbb{N}_{[1,n_c]}$ is a set-valued mapping from the parameter vector space to the index set of linear inequalities remaining in the trimmed mp-QP.

Clearly, the larger the index set $\mathbb{I}(x)$, the higher the computational cost required to solve the trimmed mp-QP of \eqref{2_6}. However, a small set $\mathbb{I}(x)$ may change the optimal solution of the original mp-QP. In fact, $\mathbb{I}(x)$ plays an important role in balancing the computational cost and optimality of the trimmed mp-QP. Intuitively, the direct design of $\mathbb{I}(x)$ in \eqref{2_6_2} requires to {\em solve} the trimmed mp-QP of \eqref{2_6} and then check its feasibility of \eqref{2_5}. If not, we should increase the size of $\mathbb{I}(x)$ and resolve the resultant mp-QP until finding the optimal solution of \eqref{2_5}.  Such a trial-and-error process is computationally demanding and not acceptable.  Instead, we adopt a general rule from \cite{nouwens2023constraint}  to certify the optimality of the trimmed mp-QP.

\subsection{On the optimality of the trimmed mp-QP} \label{sec:online_optimality}

First, we notice that $H$ is positive definite and the two constraint sets of \eqref{2_5_2} and \eqref{2_6_2} are specified by linear inequalities. This implies that the two versions of mp-QPs of \eqref{2_5} and \eqref{2_6} have unique optimal solutions, which are denoted by $z^*(x)$ and $z^*(x, \mathbb{I}(x))$, respectively. Then, we introduce an auxiliary set $\mathcal{M}(x)$ based on which an optimality condition is adopted to ensure the optimality of the trimmed mp-QP.

\begin{lemma}(\cite[Lemma 2]{nouwens2023constraint}) \label{lem:outer}
	Let the set-valued mapping $\mathbb{C}: \mathcal{R}^{n_x} \rightrightarrows \mathbb{N}_{[1,n_c]}$ denote the index set of removed inequality constraints, i.e., $\mathbb{C}(x) = \mathbb{N}_{[1,n_c]} - \mathbb{I}(x)$. If there exists a mapping $\mathcal{M}: \mathcal{R}^{n_x} \rightrightarrows \mathcal{R}^{n_z}$ such that for all $x \in \mathcal{R}^{n_x}$,
	\begin{subequations} \label{3_1}
		\begin{align}
			& z^*(x, \mathbb{I}(x)) \in \mathcal{M}(x), \label{3_1_1} \\
			& \mathcal{M}(x) \subseteq \mathcal{Z}(x, \mathbb{C}(x)), \label{3_1_2}
		\end{align}
	\end{subequations}
	the trimmed mp-QP in \eqref{2_6} does not change the optimal solution, i.e., $z^*(x, \mathbb{I}(x)) = z^*(x)$. 
\end{lemma}

\begin{figure}[t]
	\centering
	\includegraphics[width=5.6cm]{{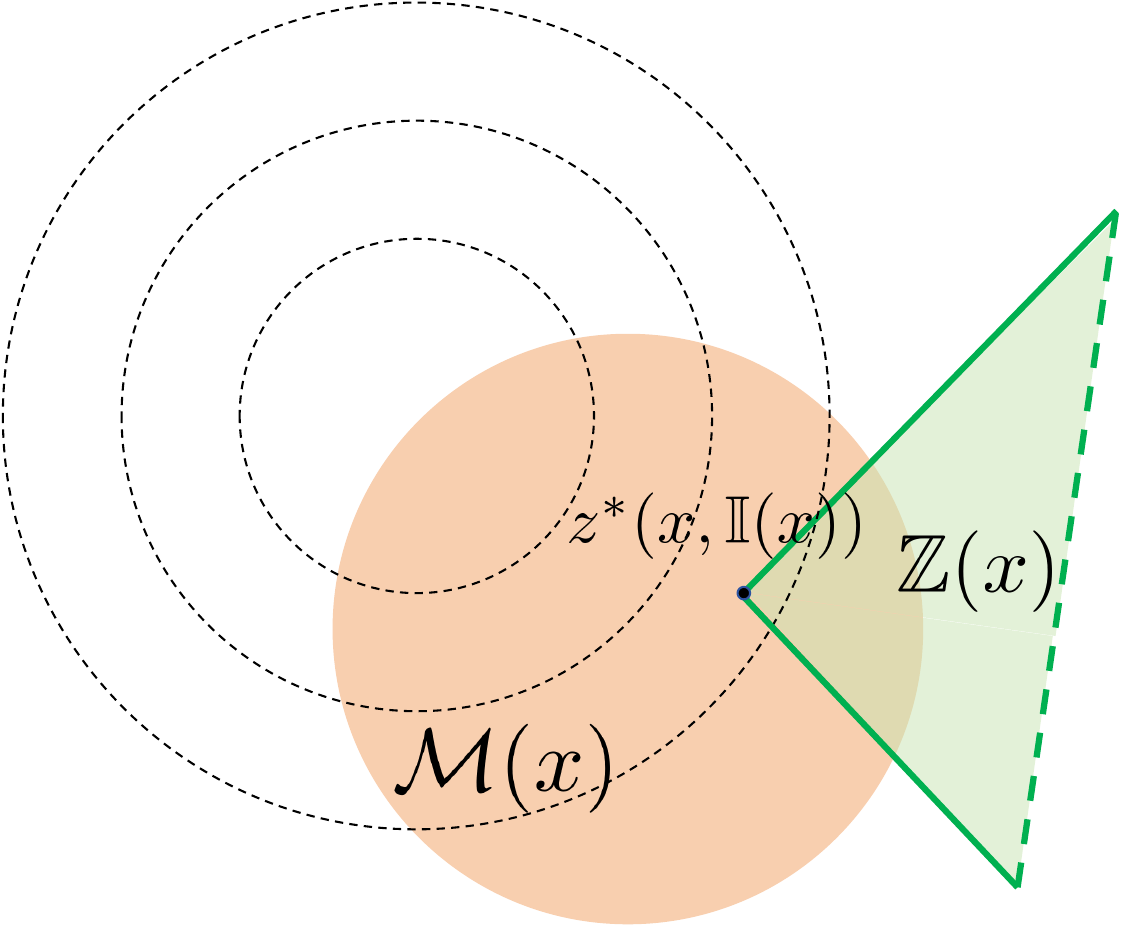}}
	\caption{An illustration of Lemma \ref{lem:outer}. The black dashed lines denote isolines of the cost function. The orange region is $\mathcal{M}(x)$, and the green solid and dashed lines represent the planes that are specified by the linear inequalities of $\mathbb{I}(x)$ and $\mathbb{C}(x)$, respectively. The left half-space of the green dashed line is the set $\mathcal{Z}(x,\mathbb{C}(x))$.}
	\label{fig:approximation}
\end{figure}

To make it self-contained, we outline the proof of Lemma \ref{lem:outer} here. In view of \eqref{3_1}, it is trivial that $z^*(x, \mathbb{I}(x)) \in \mathcal{Z}(x, \mathbb{C}(x))$.  Since $z^*(x, \mathbb{I}(x)) \in \mathcal{Z}(x, \mathbb{I}(x))$, it follows that $z^*(x, \mathbb{I}(x)) \in \mathcal{Z}(x, \mathbb{I}(x))\bigcap\mathcal{Z}(x, \mathbb{C}(x))= \mathcal{Z}(x)$. This implies that $z^*(x, \mathbb{I}(x))$ is also an optimal solution of the mp-QP in \eqref{2_5}, and thus $z^*(x, \mathbb{I}(x)) = z^*(x)$ as the optimal solution is unique.

In comparison, the joint design of a good pair of $\mathbb{I}(x)$ and $\mathcal{M}(x)$ is relatively easy, as we only need to {\em enclose} the optimal solution of \eqref{2_6} (cf. \eqref{3_1_1}) which can be achieved without solving the trimmed mp-QP.  In fact, Ref. \cite{nouwens2023constraint} leverages this idea and find an $\mathcal{M}(x)$ of \eqref{3_1} by exploiting properties of the closed-loop dynamical system in the context of linear quadratic MPC. Clearly, such a system-theoretic approach only works for the situation that the parameter vectors are sequentially generated from a stable linear time-invariant (LTI) system. More importantly, the  advantage of their method cannot be theoretically evaluated and only illustrated via simulation. In a sharp contrast, our approach of this work not only apply to {\em any} sequence of parameter vectors, but its performance  can also be explicitly evaluated. 

\subsection{Novel methods to reliably trim mp-QP from a solved instance} \label{sec:online_single}

In this subsection, we show how to utilize a solved QP($\widehat{x}$) to reliably remove redundant linear inequalities in the feasible set of \eqref{2_5}.  Specifically, we use the optimal solution $z^*(\widehat{x})$ and active set $\mathbb{A}(\widehat{x})$ of the solved QP($\widehat{x}$) where
\begin{equation}
	\mathbb{A}(\widehat{x}) := \{j \in \mathbb{N}_{[1, n_c]} \ | \ G_jz^*(\widehat{x}) = w_j +S_j\widehat{x}  \}\label{activeset}
\end{equation}
to explicitly construct a good pair of $\mathbb{I}(x)$ and $\mathcal{M}(x)$ of \eqref{3_1}.  

In view of \cite[Proposition 7.13]{rawlings2017model}, the minimizer $z^*(x)$ is Lipschitz continuous on bounded sets, implying that $z^*(x)$ can be enclosed by a ball centered at $z^*(\widehat{x})$ with a radius that is proportional to the distance  $\|\widehat{x}-x\|$ of the two parameter vectors, which is the key to the design of $\mathbb{I}(x)$ and $\mathcal{M}(x)$ of \eqref{3_1_1}.  To make it precise, we introduce the concept of global Lipschitz constant for the mp-QP. 

\begin{definition}[Global Lipschitz constant (GLC)] \label{def:Lipschitz}
	We say that $\kappa \in \mathcal{R}$ is a GLC of the mp-QP of \eqref{2_5} if for any $x_1, x_2 \in \mathcal{R}^{n_x}$ and $\forall \mathbb{I} \subseteq \mathbb{N}_{[1,n_c]}$, it holds that
	\begin{equation}
		\Vert z^*(x_1, \mathbb{I}) - z^*(x_2, \mathbb{I})  \Vert \le \kappa \Vert x_1 - x_2 \Vert. \label{3_3}
	\end{equation}
\end{definition}

\begin{remark}\label{rem_constant} Given a fixed $\mathbb{I} \subseteq \mathbb{N}_{[1,n_c]}$, it follows from \cite[Proposition 7.13]{rawlings2017model} that there exists $\kappa_{\mathbb{I}}>0$ such that \eqref{3_3} holds for $\kappa$ being replaced by $\kappa_{\mathbb{I}}$. Then, one can simply set $\kappa=\max_{\mathbb{I} \subseteq \mathbb{N}_{[1,n_c]}}\kappa_{\mathbb{I}}$, which is a uniform upper bound of Lipschitz constants over all subsets of $\mathbb{N}_{[1,n_c]}$.  \qed
\end{remark}

By \eqref{3_3}, it is trivial that $z^*(x, \mathbb{I}(x))$ lies in the ball with the center  $z^*(\widehat{x}, \mathbb{I}(x))$ and a radius $\kappa \Vert \widehat{x} - x \Vert$, i.e., 
\begin{equation}\label{ball}
	z^*(x, \mathbb{I}(x)) \in \mathcal{B}(z^*(\widehat{x}, \mathbb{I}(x)), \kappa \Vert \widehat{x} - x \Vert). 
\end{equation}

In comparison with \eqref{3_1_1}, we only need to specify the unknown ball center $z^*(\widehat{x}, \mathbb{I}(x))$ which in fact is exactly the optimal solution of the solved QP($\widehat{x}$), i.e., $z^*(\widehat{x})=z^*(\widehat{x}, \mathbb{I}(x))$, if the index set $\mathbb{I}(x)$ is selected to satisfy $ \mathbb{A}(\widehat{x})\subseteq \mathbb{I}(x)$. Formally, we have the following result. 

\begin{lemma} \label{lem:equal}
	Given two parameter vectors $x, \widehat{x}\in \mathcal{R}^{n_x}$, if $\mathbb{A}(\widehat{x}) \subseteq \mathbb{I}(x) \subseteq  \mathbb{N}_{[1,n_c]}$, then $z^*(\widehat{x}, \mathbb{I}(x)) = z^*(\widehat{x})$. 	
\end{lemma}
\begin{pf}
	Since $\mathbb{A}(\widehat{x}) \subseteq \mathbb{I}(x) \subseteq \mathbb{N}_{[1,n_c]}$, it holds that
	\begin{equation}
		V(z^*(\widehat{x}, \mathbb{A}(\widehat{x}))) \le V(z^*(\widehat{x}, \mathbb{I}(x))) \le V(z^*(\widehat{x})). \label{3_1_P0}
	\end{equation}
	
	By the definition of active set in \eqref{activeset}, it follows from \cite{rawlings2017model}  that $V(z^*(\widehat{x}, \mathbb{A}(\widehat{x}))) = V(z^*(\widehat{x}))$. Jointly with \eqref{3_1_P0}, it holds
	\begin{equation}
		V(z^*(\widehat{x}, \mathbb{A}(\widehat{x}))) =  V(z^*(\widehat{x}, \mathbb{I}(x))) = V(z^*(\widehat{x})). \label{3_1_P1}
	\end{equation}
	
	Since both $z^*(\widehat{x}, \mathbb{I}(x))$ and $z^*(\widehat{x})$ belong to $\mathcal{Z}(\widehat{x}, \mathbb{A}(\widehat{x}))$, they are two optimal solutions of the trimmed instance of \eqref{2_6} with $\mathbb{I}(x) = \mathbb{A}(\widehat{x})$.  By the uniqueness of its optimal solution, it follows that
	\begin{equation}
		z^*(\widehat{x}, \mathbb{I}(x)) = z^*(\widehat{x})\label{3_1_P2} 
	\end{equation}
	which completes the  proof. 
\end{pf}

In view of \eqref{ball} and Lemma \ref{lem:equal}, let $\mathcal{M}(x) = \mathcal{B}(z^*(\widehat{x}), \kappa \Vert x - \widehat{x} \Vert)$. Then, \eqref{3_1_1} is insured as
\begin{equation}
	z^*(x, \mathbb{I}(x)) \in\mathcal{M}(x)~\text{if}~\mathbb{A}(\widehat{x})\subseteq \mathbb{I}(x).  \label{acondition}
\end{equation}
Note that it is easy to construct the pair of $\mathbb{I}(x)$ and $\mathcal{M}(x)$ as both $z^*(\widehat{x})$ and $\mathbb{A}(\widehat{x})$ are directly obtained from the solved QP($\widehat{x}$).  

Next, we focus on the condition of \eqref{3_1_2}  by further increasing the index set $\mathbb{I}(x)$.  To elaborate it, we propose the following index set
\begin{equation} \label{3_8}
	\mathbb{D}(x) = \{j \in \mathbb{N}_{[1, n_c]} \ |  \mathcal{M}(x) \subseteq \mathcal{Z}_j(x) \}.
\end{equation}

Since $\mathcal{M}(x)$ is a ball and $\mathcal{Z}_j(x)$ is a half-plane, then $\mathbb{D}(x)$ is easy to obtain as well and will be detailed in Remark \ref{rem_formula}. Obviously, \eqref{3_8} implies that $\mathcal{M}(x)\subseteq \mathcal{Z}(x,\mathbb{D}(x))$ and \eqref{3_1_2} is insured as
\begin{equation}
	\mathcal{M}(x)\subseteq \mathcal{Z}(x,\mathbb{C}(x))~\text{if}~ \mathbb{C}(x)\subseteq \mathbb{D}(x).  \label{bcondition}
\end{equation}

Jointly with $\mathbb{C}(x) = \mathbb{N}_{[1,n_c]} - \mathbb{I}(x)$, the condition \eqref{3_1} holds for the pair of sets  $\mathcal{M}(x) = \mathcal{B}(z^*(\widehat{x}), \kappa \Vert x - \widehat{x} \Vert)$ and 
\begin{equation}\label{indexseti}
	\begin{aligned}
		\mathbb{I}(x)&= \mathbb{A}(\widehat{x})\bigcup \big(\mathbb{N}_{[1, n_c]}-\mathbb{D}(x)\big)\\
		&=\mathbb{A}(\widehat{x})\bigcup \big(\mathbb{A}^c(\widehat{x})-\mathbb{D}(x)\big)
	\end{aligned}
\end{equation} 
where the so-called inactive set $\mathbb{A}^c(\widehat{x})$ is the complement of $\mathbb{A}(\widehat{x})$ and given as 
$$\mathbb{A}^c(\widehat{x}):= \{j \in \mathbb{N}_{[1, n_c]} \ | \ G_jz^*(\widehat{x}) < w_j +S_j\widehat{x} \}.$$

\begin{remark}\label{rem_activeset} If the parameter vector $x$ is  close to $\widehat{x}$, this implies that QP($x$) is also  ``close'' to the solved QP($\widehat{x}$) and it is expected that the results of QP($\widehat{x}$) could provide useful information to trim QP($x$).  This can be justified from our design of $\mathbb{I}(x)$ and $\mathcal{M}(x)$. Specifically, if $\Vert x - \widehat{x} \Vert$ is small, the set $\mathcal{M}(x)$ tends to be a small ball centered at $z^*(\widehat{x})$. Jointly with that $z^*(\widehat{x})$ is an interior point of $\mathcal{Z}_j(x)$ for any $j\in  \mathbb{A}^c(\widehat{x})$, then  $\mathcal{M}(x) \subseteq \mathcal{Z}_j(x)$ could be potentially satisfied. This indicates that $\mathbb{D}(x)\approx  \mathbb{A}^c(\widehat{x})$ and thus $\mathbb{I}(x)\approx \mathbb{A}(\widehat{x})$ (c.f. \eqref{indexseti}). In this case, the feasible set in the trimmed QP($x$) of \eqref{2_6}  tends to be  $\mathcal{Z}(x, \mathbb{A}(\widehat{x}))$, which is the minimum set of linear inequalities required to retain the optimal solution at $x=\widehat{x}$. \end{remark}

\begin{algorithm}[t]
	\caption{Reliably learn to trim new mp-QP instance of \eqref{2_5} from a solved one} \label{alg:brief}
	\begin{algorithmic}[1]
		\Require a GLC $\kappa$ of \eqref{2_5}, and a pair $\left(z^*(\widehat{x}),\mathbb{A}(\widehat{x})\right)$ from the solved QP($\widehat{x}$)
		\Ensure a trimmed QP in the form of \eqref{2_6}
		\State Compute \eqref{3_11};  \label{alg_1_1} 
		\State Set $\mathbb{I}(x) =\mathbb{A}(\widehat{x})\bigcup\left( \mathbb{A}^c(\widehat{x})-\mathbb{D}(x)\right)$. \label{alg_1_3}
	\end{algorithmic}
\end{algorithm}
\begin{figure}[t]
		\centering
	\includegraphics[width=9cm]{{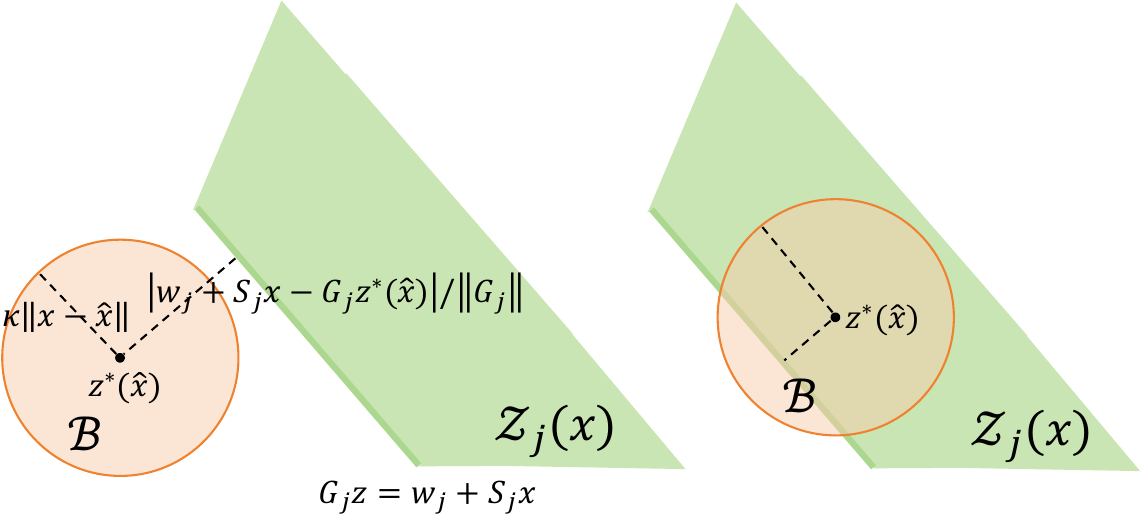}}
		\caption{A geometrical illustration of $\mathcal{B}(z^*(\widehat{x}), \kappa \Vert x - \widehat{x} \Vert) \not\subseteq \mathcal{Z}_j(x)$. }
		\label{fig:ball}
	\end{figure}
\begin{remark}\label{rem_formula} 
	We provide an explicit formula to compute  $\mathbb{D}(x)$ in \eqref{3_8}. It is trivial that $\mathcal{B}(z^*(\widehat{x}), \kappa \Vert x - \widehat{x} \Vert) \subseteq \mathcal{Z}_j(x)$ of \eqref{3_8}  {\em if and only if} (a) $z^*(\widehat{x})\in \mathcal{Z}_j(x)$, i.e., $G_j z^*(\widehat{x}) \le w_j + S_j x$, and (b) the distance from the ball center $z^*(\widehat{x})$ to the hyperplane of $\mathcal{Z}_j(x)$ should be greater than the ball radius $ \kappa \Vert x - \widehat{x} \Vert$, i.e., 
	$
	|w_j + S_j x - G_j z^*(\widehat{x})|/\Vert G_j \Vert  \ge \kappa \Vert x - \widehat{x} \Vert.
	$ See Fig. \ref{fig:ball} for an illustration of $\mathcal{B}(z^*(\widehat{x}), \kappa \Vert x - \widehat{x} \Vert) \not\subseteq \mathcal{Z}_j(x)$. Then, the set $\mathbb{A}^c(\widehat{x})-\mathbb{D}(x)$ of \eqref{indexseti} can be explicitly computed as
	%
	%
	%
	\begin{equation} \label{3_11}
		\begin{aligned} 
			&\mathbb{A}^c(\widehat{x})-\mathbb{D}(x) \\ 
			&=\left \{j \in \mathbb{A}^c(\widehat{x}) \ | 
			 \kappa  \cdot \Vert x - \widehat{x} \Vert > \frac{w_j + S_j x - G_j z^*(\widehat{x})}{\Vert G_j \Vert }\right\}.
		\end{aligned}
	\end{equation}

\end{remark}

\begin{remark}\label{rem_glc}
	In light of  \eqref{3_11}, the smaller the GLC $\kappa$, the smaller the number of  elements of $\mathbb{A}^c(\widehat{x})-\mathbb{D}(x)$ tends to be and thus the slimmer the trimmed mp-QP of \eqref{2_6}. In Section \ref{sec:offline_lipschitz}, we shall provide a method to compute GLC. 
\end{remark}

Before closing this subsection, we summarize our idea to trim QP($x$) in Algorithm \ref{alg:brief}. By Lemma \ref{lem:outer}, it is straightforward that the trimmed mp-QP of Algorithm \ref{alg:brief} does not change the optimal solution, which is formalized below. 
\begin{theorem} \label{thm:data}
	The trimmed mp-QP of Algorithm \ref{alg:brief} has the same optimal solution as that of the original mp-QP of \eqref{2_5}, i.e., $z^*(x, \mathbb{I}(x)) = z^*(x)$.
\end{theorem}

\subsection{Learn to trim from multiple solved mp-QP instances} \label{sec:online_multiple}
Now, we use a sequence of pairs $\{(z^*(\widehat{x}^k), \mathbb{A}(\widehat{x}^k))\}_{k\in  \mathbb{N}_{[1,q]}}$ from multiple solved mp-QP instances to refine $\mathbb{I}(x)$ and $\mathcal{M}(x)$ of \eqref{3_1}.  Naively, one may use the results of the instance with the closest parameter vector via Algorithm \ref{alg:brief}, i.e., $\widehat{x}\in\arg\min_{\widehat{x}^k}\{\|x-\widehat{x}^k\|, k\in  \mathbb{N}_{[1,q]}\}$. Though simple, such an idea does not fully exploit the results of all solved instances. 

A straightforward approach is to apply Algorithm \ref{alg:brief} to each solved QP($\widehat{x}^k$) in parallel and then take the intersection of their resultant index sets $\mathbb{I}^k(x)$, but this may be problematic as illustrated below.

\begin{example} \label{exa:mul}
	Consider a one-dimensional mp-QP  with $V(z)= z^2 + x\cdot z,   \mathcal{Z}(x)=\{z| \ z \le x, z \le -x - 4 \}$. Clearly, $\kappa = 1$ is a GLC of this problem.  Suppose that we have solved $\text{QP}(-1)$ and $\text{QP}(-3)$, and obtain that
	\begin{equation} \notag 
		\begin{aligned}
			& z^*(-1) = -3, \mathbb{A}(-1) = \{2\}~\text{and}~ \\
			& z^*(-3) = -3, \mathbb{A}(-3) = \{1\}.
		\end{aligned}
	\end{equation}
	
	Then, we use them to trim $\text{QP}(-2)$ via Algorithm \ref{alg:brief},  respectively, and obtain two trimmed instances with index sets of linear inequalities as $\mathbb{I}^1(-2) = \{2\}$ and $\mathbb{I}^2(-2) = \{1\}$. One can easily verify that both trimmed instances preserve the optimal solution, i.e., $z^*(-2, \mathbb{I}^1(-2)) = z^*(-2, \mathbb{I}^2(-2)) = z^*(-2)$. If we simply take an intersection of the two index sets, the trimmed instance with $ \mathbb{I}(-2)= \mathbb{I}^1(-2) \cap \mathbb{I}^2(-2) = \emptyset$ cannot preserve the optimal solution of $\text{QP}(-2)$, i.e., $z^*(-2, \mathbb{I}(-2))\neq z^*(-2)$. \qed
\end{example}

To remedy the issue in Example \ref{exa:mul}, we require the following linear independence constraint qualification (LICQ) condition.
\begin{assumption} \label{ass:LICQ}
	(LICQ  \cite{borrelli2017predictive}) For any parameter vector $x \in \mathcal{R}^{n_x}$, the rows of $G_{\mathbb{A}(x)}$ are linearly independent.
\end{assumption}

Under the LICQ condition of Assumption \ref{ass:LICQ}, the optimal solution can be specified by a unique minimum set of linear inequalities. That is, the trimmed mp-QP of \eqref{2_6} with $\mathbb{I}(x)=\mathbb{A}(x)$ is the slimmest one that preserves the optimal solution of \eqref{2_5}. This minimum set of linear inequalities must be included in any trimmed mp-QP. Conversely, any trimmed mp-QP containing this set retains the same optimal solution as established in the proof of Lemma \ref{lem:equal}. Thus, the intersection of all index sets of $\mathbb{I}^k(x)$ obtained from Algorithm \ref{alg:brief} over the solved QP($\widehat{x}^k$), will preserve the optimal solution. Importantly, Assumption \ref{ass:LICQ} guarantees that the cardinality of the active set does not exceed $n_z$, i.e., $|\mathbb{A}({x})| \le n_z$. When LICQ is not satisfied, we propose a practical remedy: select the instance corresponding to the closest parameter vector, as implemented in Section~\ref{sec:mpc_offline} (see, e.g., \eqref{mindist}). We formalize the complete procedure in Algorithm \ref{alg:brief_mul} and  establish its optimality in the following proposition. 

\begin{theorem} \label{pro:mul} Under Assumption \ref{ass:LICQ}, the trimmed mp-QP of Algorithm \ref{alg:brief_mul} has the same optimal solution as that of the original mp-QP of \eqref{2_5}, i.e., $z^*(x, \mathbb{I}(x)) = z^*(x)$. 
\end{theorem}
\begin{pf}
	See Appendix \ref{app:mul}.
\end{pf}
\begin{algorithm}[t]
	\caption{Reliably learn from multiple solved mp-QP instances to trim new instance} \label{alg:brief_mul}
	\begin{algorithmic}[1]
		\Require a GLC $\kappa$ of \eqref{2_5}, and a sequence of pairs $\left\{\left(z^*(\widehat{x}^k), \mathbb{A}(\widehat{x}^k)\right)\right\}_{k\in  \mathbb{N}_{[1,q]}}$ from solved QP($\widehat{x}^k$)
		\Ensure a trimmed QP in the form of \eqref{2_6}
		\State Set $k = 1$ and $\mathbb{I}(x) = \mathbb{N}_{[1, n_c]}$.
		\While{$k \le q$}
		\State \hspace{-.65cm} Compute 
		\begin{align*} 
			\mathbb{T}(x,\widehat{x}^k) &=\left \{j \in \mathbb{A}^c(\widehat{x}^k)\bigcap \mathbb{I}(x) \ | \right. \\ 
			&\left.~~~\kappa  \cdot \Vert x - \widehat{x}^k \Vert > \frac{w_j + S_j x - G_j z^*(\widehat{x}^k) }{\Vert G_j \Vert}\right\}.
		\end{align*} 
		\State Set $\mathbb{I}(x) \leftarrow\left(\mathbb{A}(\widehat{x}^k)\bigcap \mathbb{I}(x)\right)\bigcup \mathbb{T}(x,\widehat{x}^k)$. \label{alg_2_3}
		\State $k \leftarrow k + 1$.
		\EndWhile
	\end{algorithmic}
\end{algorithm}

Similar to Algorithm 1, each loop of Algorithm 2 uses a solved QP($\widehat{x}^k$) to identify redundant constraints. The main difference from Algorithm 1 is that the follow-up loop only needs to focus on the remaining set of inequality constraints $\mathbb{I}(x)$ rather than the whole set of constraints $\mathbb{N}_{[1, n_c]}$. In fact,  Algorithm \ref{alg:brief_mul} can also be implemented in a parallel way where each processor $k$ individually computes an index set $\mathbb{I}^k(x)$ as 
\begin{subequations}
	\begin{align}
		\mathbb{R}(x,\widehat{x}^k) &:= \left\{j \in \mathbb{A}^c(\widehat{x}^k) \ | \right. \notag \\ 
		&\left.~~~\kappa  \cdot \Vert x - \widehat{x}^k \Vert > \frac{w_j + S_j x - G_j z^*(\widehat{x}^k) }{\Vert G_j \Vert}\right\}, \label{3_12_1} \notag \\
		\mathbb{I}^k(x) &= \mathbb{A}(\widehat{x}^k)\bigcup \mathbb{R}(x,\widehat{x}^k)\notag
	\end{align}
\end{subequations}
and then aggregates as $\mathbb{I}(x)=\bigcap_{k=1}^q \mathbb{I}^k(x)$. 

\section{Performance Analysis of Algorithm \ref{alg:brief}} \label{sec:offline}
In this section, we provide a method to compute the GLC via MILP. Then, we quantify the number of linear inequalities in the trimmed mp-QP instance of Algorithm \ref{alg:brief}.


\subsection{Computation method of global Lipschitz constant} \label{sec:offline_lipschitz}
In Remark \ref{rem_glc}, we have briefly shown that the GLC is essential to the efficiency of our trimming Algorithm \ref{alg:brief}. However, finding the minimum one  is NP-hard \cite{9123670, d54e99a6, 9928332, teichrib2023efficient}. In this subsection, we provide a certifiable GLC via mixed integer linear program (MILP).  To this goal, we define the minimum GLC under $l_a$-norm  as $\kappa_a$, $a \in \{1, 2, \infty\}$, i.e., $\kappa_a$ is obtained by replacing $l_2$-norm in \eqref{3_3} with $l_a$-norm.

\begin{lemma} \label{lem:LipschitzPro} The minimum GLCs satisfy that 
 \begin{equation} \notag
	 	\max \left\{  \frac{1}{\sqrt{n_z}} \kappa_1, \frac{1}{\sqrt{n_x}} \kappa_\infty  \right\} \le\kappa_2 \le  \sqrt{\kappa_1\kappa_\infty}
	 \end{equation}
	 and a GLC of \eqref{3_3} can be given as 
	\begin{equation} \label{eq:kappa}
		\kappa = \sqrt{\kappa_1\kappa_\infty}.
	\end{equation}
\end{lemma}

\begin{pf}
	Since the Hessian matrix $H$ is positive definite, the trimmed instance of \eqref{2_6} with $\mathbb{I}(x)=\mathbb{I}$ admits a piecewise affine function $z^*(x, \mathbb{I})$ over $x$ \cite{bemporad2002explicit}. We define $K(x,\mathbb{I})$ as its slope at the parameter vector $x$.
	First, we derive an upper bound of $\kappa_2$ as follows:
	\begin{equation} \notag 
		\begin{aligned}
			\kappa_2=\max_{x, \mathbb{I}}\Vert K(x,\mathbb{I}) \Vert^2 & \le \max_{x, \mathbb{I}}\Vert K(x,\mathbb{I}) \Vert_1 \Vert K(x,\mathbb{I}) \Vert_{\infty} \\
			& \le \max_{x, \mathbb{I}}\Vert K(x,\mathbb{I}) \Vert_1 \max_{x, \mathbb{I}} \Vert K(x,\mathbb{I}) \Vert_{\infty} \\
			& \le \kappa_1 \kappa_\infty
		\end{aligned}
	\end{equation}
     where the first equality follows from \cite[Section 2.1]{darup2017maximal}.
	Then, we establish a lower bound of $\kappa_2$. From the equivalence of matrix norms \cite[Chapter 5]{horn2012matrix},  i.e.,
	\begin{equation} \notag
		\Vert  K(x,\mathbb{I})\Vert \ge \max \left\{ \frac{1}{\sqrt{n_z}} \Vert  K(x,\mathbb{I}) \Vert_1, \frac{1}{\sqrt{n_x}} \Vert  K(x,\mathbb{I}) \Vert_\infty \right\},
	\end{equation}
	 the minimum GLC $\kappa_2$ can be lower bounded by
	\begin{equation} \notag
		\kappa_2 \ge \max \left\{  \frac{1}{\sqrt{n_z}} \kappa_1, \frac{1}{\sqrt{n_x}} \kappa_\infty  \right\}.
	\end{equation}
	Finally, we choose the upper bound as a GLC, i.e., $\kappa = \sqrt{\kappa_1\kappa_\infty}$.
\end{pf}

We next provide a method to compute $\kappa_1$ and $\kappa_\infty$ in \eqref{eq:kappa}.
\begin{lemma} \label{lem:Lipschitz}
	Both $\kappa_1$ and $\kappa_\infty$ can be expressed  via the solutions of MILPs.
\end{lemma}
\begin{pf}
	See Appendix \ref{app:Lipschitz}.
\end{pf}

\subsection{Performance evaluation of Algorithm \ref{alg:brief}}  \label{sec:offline_properties}
In this subsection,  we evaluate the performance of Algorithm \ref{alg:brief} by quantifying the number of linear inequalities in the trimmed instance of \eqref{2_6}.  As noted in Remark \ref{rem_activeset}, the closeness of $\widehat{x}$ to $x$ is essential to the efficiency of removing redundant inequalities in QP($x$).  Thus, we adopt an increasing sequence of $\{\sigma_i\}$, which are specified as the infimum values of MILP problems in the sequel, to classify the closeness between two parameter vectors, and use them to evaluate the quality of QP($\widehat{x}$) in trimming QP($x$). To this end, it follows from \eqref{indexseti} that the key is to quantify the cardinality of $\mathbb{D}(x)$. Since it depends on both $x$ and $\widehat{x}$, i.e., 
\begin{equation} \label{Dx}
	\mathbb{D}(x) = \left\{j \in \mathbb{N}_{[1, n_c]} \ | \mathcal{B}(z^*(\widehat{x}), \kappa\Vert x - \widehat{x} \Vert) \subseteq \mathcal{Z}_j(x)\right\},
\end{equation}
we do not directly work on it, and instead lift both $\mathcal{Z}_j(x)$ and $\mathcal{B}(z^*(\widehat{x}), \kappa\Vert x - \widehat{x} \Vert)$ to characterize their relation, based on which we provide a lower bound for the cardinality of $\mathbb{D}(x)$. 

Consider a lifted version of the feasible set as
\begin{equation}
	\mathcal{V} = \{v \ | \ Dv \le w\},
\end{equation}
where $v=[x',z']'\in\mathcal{R}^{n_v}$, $D = [-S, G] \in \mathcal{R}^{n_c \times n_v}$ and $n_v = n_x + n_z$. Accordingly, let $\mathcal{V}_j = \{v \ | \ D_j v \le w_j \}$ for any $j \in \mathbb{N}_{[1, n_c]}$ and define
\begin{subequations}\label{4_10sec}
	\begin{align}
		v(\widehat{x})&= \begin{bmatrix} \widehat{x}^T,  z^*(\widehat{x})^T \end{bmatrix}^T \in \mathcal{R}^{n_v}, \label{4_10_1} \\
		\mathbb{E}(x) &= \left\{j \in \mathbb{N}_{[1, n_c]} \ |\  \mathcal{B}(v(\widehat{x}), \sqrt{1 + \kappa^2}\Vert x - \widehat{x} \Vert) \subseteq \mathcal{V}_j \right \}. \label{4_10_2} 
	\end{align}
\end{subequations}

Note that $\mathcal{V}=\bigcap\nolimits_{j \in \mathbb{N}_{[1, n_c]}}\mathcal{V}_j$ is also a non-empty polyhedron, and we obtain following results.
\begin{lemma} \label{lem:extended} The sets in \eqref{4_10sec} satisfy the following relation
	\begin{enumerate} [(a)] \label{4_10}
		\item \label{4_10_3} $v(\widehat{x}) \in \mathcal{V}$.
		\item \label{4_10_4} $\mathbb{E}(x) \subseteq \mathbb{D}(x)$.
	\end{enumerate}
\end{lemma}

\begin{pf}
	See Appendix \ref{app:extended}.
\end{pf}	
By Lemma \ref{lem:extended}, it immediately implies that $v(\widehat{x})\in \mathcal{V}_j$, i.e., 
the ball center $v(\widehat{x})$ of \eqref{4_10_2} automatically lies in the lifted half-space $\mathcal{V}_j$ for any $j \in \mathbb{N}_{[1, n_c]}$, which is not the case for the set $\mathbb{D}(x)$ of \eqref{Dx}. From this view, it  seems easier to numerate the element of $\mathbb{E}(x)$ than that of $\mathbb{D}(x)$.  Since $v(\widehat{x}) \in \mathcal{V}$ holds automatically, we drop the dependence on the particular ball center in $\mathbb{E}(x)$ by noting the following relation  
$$
|\mathbb{E}(x)|\ge \inf_{v\in \mathcal{V}} \left | \left\{j \in \mathbb{N}_{[1, n_c]} \ |\  \mathcal{B}(v, \sqrt{1 + \kappa^2}\cdot \Vert x - \widehat{x} \Vert) \subseteq \mathcal{V}_j\right \}\right |
$$
and then solely  focus on the effect of ball radius via studying
\begin{equation}
	\inf_{v\in \mathcal{V}} \left | \left\{ j \in \mathbb{N}_{[1, n_c]} \ |\  \mathcal{B}(v, r) \subseteq \mathcal{V}_j \right \}\right | \in  \mathbb{N}_{[0,n_c]}, \forall r\ge 0.
\end{equation}

Now, we are ready to specify the sequence of $\{\sigma_i\}$ to classify the distance between two parameter vectors, and then use them to find a lower bound on the cardinality of $\mathbb{E}(x)$.

\begin{lemma} \label{lem:D}
	For any $i \in \mathbb{N}_{[1,n_c]}$, let $\sigma_i$ be the maximum value of the following optimization problem 	
	\begin{subequations} \label{4_11}
		\begin{align}
			\hspace{-0.2cm} \sigma_i= & \max_{r\ge 0} \  r,  \label{4_11_1} \\
			& \text{s.t.}~\min_{v \in \mathcal{V}} \left | \left\{ j \in \mathbb{N}_{[1, n_c]} \ | \ \mathcal{B}(v, r) \subseteq \mathcal{V}_j \right \}\right | \ge n_c-i. \label{4_11_2} 
		\end{align}
	\end{subequations}
	Then, $\sigma_i$ is finite and satisfies that
	\begin{enumerate}[(a)]
		\item \label{4_11_3} $0 \le \sigma_{i} \le \sigma_{i+1}$.
		\item \label{4_11_5} If $v \in \mathcal{V}$ and $0 \le r \le \sigma_{i}$, then
		\begin{equation}
			\left | \left\{ j\in \mathbb{N}_{[1, n_c]} \ | \ \mathcal{B}(v, r) \subseteq \mathcal{V}_j \right \}\right |  \ge n_c-i. \label{4_11_6}
		\end{equation}
	\end{enumerate}
\end{lemma}

\begin{pf}
	See Appendix \ref{app:D}.
\end{pf}	

Note that the attainability of maximization and minimization operators in \eqref{4_11} are proved in Appendix \ref{app:D}. The problem of \eqref{4_11}   does not depend on  parameter vectors, and  can be computed in prior.  Moreover, the optimal value of \eqref{4_11}  can be obtained via solving an MILP, which is formalized below. 

\begin{lemma} \label{lem:MILP}
	The optimal value of \eqref{4_11} can be obtained via solving the following  MILP. 
	\begin{equation} \label{4_12}
		\begin{aligned}
			\sigma_i= &\inf_{v \in \mathcal{V}, r\ge 0} \  r, \\
			&\text{s.t.}~ 
			\frac{w_j - D_j v}{\Vert D_j \Vert} -M  \delta_j < r \le \frac{w_j - D_j v}{\Vert D_j \Vert} + M (1 - \delta_j), \\
			& \sum_{j=1}^{n_c} \delta_j< n_c - i, \delta_j \in \{0,1\}, \ j \in \mathbb{N}_{[1, n_c]}.
		\end{aligned}
	\end{equation}
	where $M$ is a big number \cite{9928332}. 
\end{lemma}
\begin{pf}See Appendix \ref{app:MILP}. 
\end{pf}

Though the infimum in \eqref{4_12} may not be attainable, it does not incur any trouble as \eqref{4_11_6} only uses the value of $\sigma_i$. 

Now, we use Lemma \ref{lem:D} to quantify the number of linear inequalities in the trimmed mp-QP instance of \eqref{2_6}. 
\begin{theorem} \label{thm:distance}
	Suppose that Assumption \ref{ass:LICQ} holds.  If $\Vert x - \widehat{x}\Vert \le {\sigma_i}/{\sqrt{1 + \kappa^2}}$ for some $i \in \mathbb{N}_{[1,n_c]}$, then $|\mathbb{I}(x)| \le n_z + i $. 
\end{theorem}

\begin{pf}
	By Lemma \ref{lem:extended}(\ref{4_10_4}) and Lemma \ref{lem:D}(\ref{4_11_5}), it holds that  $|\mathbb{D}(x)| \ge |\mathbb{E}(x)| \ge n_c - i $.
	In addition, Assumption \ref{ass:LICQ} implies that $|\mathbb{A}(\widehat{x})| \le n_z$.  Then, it follows from \eqref{indexseti} that
	\begin{equation}  \notag
		|\mathbb{I}(x)| \le n_z + n_c - |\mathbb{D}(x)| \le n_z + i 
	\end{equation}  
	which completes the proof. 
\end{pf}

Theorem \ref{thm:distance} shows that if the parameter vector $\widehat{x}$ of the solved  $\text{QP}(\widehat{x})$ satisfies $\|x-\widehat{x}\|\le \sigma_i/\sqrt{1+\kappa^2}$, then the number of linear inequalities  in the trimmed instance of Algorithm \ref{alg:brief} can be made at most $n_z + i$. It should be noted that this result is generically conservative. If the number of linear inequalities $n_c$ in QP($x$) is significantly greater than the dimension of decision vector $n_z$, which is usually the case in many applications \cite{britzelmeier2020nonsmooth, wu2021coordinated}, the number of removed inequalities ($n_c-|\mathbb{I}(x)|$) is significant and thus trimming mp-QP is very rewarding. 
 
\section{Application to the linear quadratic MPC} \label{sec:bridge}
In this section, we shall learn to trim mp-QP of the linear quadratic MPC whose parameter vector is the encountered state vector of a controlled system. Firstly, we introduce the linear MPC in the form of mp-QP. Then, we adaptively trim mp-QP instances in the closed-loop system of MPC and show that the number of linear inequalities  in the trimmed QP($x_k$) decrease to zero after a finite timestep. Finally, we utilize offline solved instances to trim MPC.

\subsection{The linear quadratic MPC}

We consider the stabilization problem of a constrained linear time-invariant (LTI) system
\begin{equation} \label{2_1}
	x_{k+1} = Ax_k + Bu_k
\end{equation}
with the state vector $x_k \in \mathcal{X}\subseteq \mathcal{R}^{n}$, the control input vector $u_k \in \mathcal{U}\subseteq \mathcal{R}^{m}$. The MPC   is a well known methodology for synthesizing control laws of \eqref{2_1} under the prescribed operating constraints, where the control action is obtained via solving a finite horizon open-loop optimal control problem \cite[Chapter 2]{rawlings2017model}, i.e.,  
\begin{equation}\label{2_2}
	\begin{aligned}
			\mathop{\text{min.}}\limits_{u_{\cdot|k}} & \ J(x_{.|k}, u_{.|k}) \\
			\text{s.t.}~ & x_{t+1|k} = Ax_{t|k} + Bu_{t|k}, \\
			& x_{t|k} \in \mathcal{X}, u_{t|k} \in \mathcal{U}, t \in \mathbb{N}_{[0, N-1]}, \\
			& x_{N|k} \in \mathcal{X}_N, \\
			& x_{0|k} = x_k.
		\end{aligned}
\end{equation}

In this work,  we assume that  the above $\mathcal{X}$ and $\mathcal{U}$, and $\mathcal{X}_N$ are bounded polyhedral sets.  Moreover,  the objective function  takes the following quadratic form
\begin{equation}
	J(x_{.|k}, u_{.|k}) = \sum^{N-1}_{t=0} (x_{t|k}^TQx_{t|k} + u_{t|k}^TRu_{t|k}) + x_{N|k}^TPx_{N|k}, \label{2_3}
\end{equation} with positive definite matrices $Q, R, P$.  Since the predicted state vectors $x_{.|k}$ can be linearly expressed by $u_{.|k}$ and $x_k$, the linear quadratic MPC in \eqref{2_2} can be reformulated as a mp-QP in the form of \eqref{2_5}   parameterized by $x_k$ \cite{bemporad2002explicit}. 

Starting from an initial state $x_0 \in \mathcal{X}$,  the MPC only applies the first input vector of the optimal sequence of \eqref{2_2} to the system \eqref{2_1}, i.e., 
\begin{equation}
	u_{\text{MPC}}(x_k) = u^*_{0|k}. \label{2_4}
\end{equation}

Under standard design choices, see e.g., \cite[Chapter 2]{rawlings2017model}, the closed-loop system exponentially converges to the origin.  In this work, we do not include the design details and simply make the following assumption.

\begin{assumption} \label{ass:basic}
	For the LTI system \eqref{2_1} under the control law \eqref{2_4}, it holds that
	\begin{enumerate}[(a)]
			\item The mp-QP in \eqref{2_2} is recursively feasible starting from any $x_0 \in \mathcal{X}$.
			\item The closed-loop system is exponentially stable, i.e., there exist positive constants $c>0$ and $\beta\in (0,1)$ such that $\Vert x_k \Vert \le c\Vert x_0 \Vert \beta^k$.
		\end{enumerate}
\end{assumption}

For simplicity, we abuse notation in this section and simply assume that the mp-QP of MPC in \eqref{2_2} is directly expressed in the form of \eqref{2_5} where  $n_z=Nm$, $n_x=n$, $n_c = \vert \mathcal{X}_N \vert +N(\vert \mathcal{X} \vert+\vert \mathcal{U} \vert)$ under the condensed approach \cite{9928332} and $\vert \mathcal{X} \vert$ denotes the number of linear inequalities to specify the polyhedral set $\mathcal{X}$. The goal of this section is to trim the QP($x_k$) of \eqref{2_5} where the parameter vector $x_k$ is the state vector of the constrained LTI system \eqref{2_1} under the control law \eqref{2_4}.

Different from the mp-QP of Section \ref{sec:online}, we are interested in trimming mp-QP instances in the closed-loop system whose state vector decreases to the origin (cf. Assumption \ref{ass:basic}). This suggests that both the state difference $\|x_k-x_{k-1}\|$ and the number of linear inequalities in the trimmed QP($x_k$) are expected to decrease to zero. This section rigorously provides an affirmative answer and shows that eventually we only need to solve a {\em constraint-free} instance.

\subsection{Adaptively trim mp-QP instances in the closed-loop system of MPC }\label{subsec_adaptive}
In this subsection, we use our trimming methods in Section \ref{sec:online} to adaptively remove redundant inequality constraints in the QP($x_k$) of MPC, which is summarized in Algorithm \ref{alg:adaptive}. Then, we evaluate its performance by quantifying the number of linear inequalities in its trimmed version. 


It should be noted that Algorithm \ref{alg:adaptive}  requires to store the results of instance at the last timestep.  As the closed-loop system exponentially converges to the origin (cf. Assumption \ref{ass:basic}), it is expected that the number of inequality constraints decreases and even to zero if $0\in\text{int}(\mathcal{X}), 0\in\text{int}(\mathcal{U})$ and $0\in\text{int}(\mathcal{X}_N)$. Informally, if the origin lies in the boundary of a constraint set, there is at least one active linear inequality for the optimal solution of QP($0$) of MPC, and such an active inequality usually cannot be removed. See the discussion after Assumption \ref{ass:LICQ}. 

\begin{algorithm}[!t]
	\caption{Adaptively trim MPC in the closed-loop system} \label{alg:adaptive}
	\begin{itemize}
		\item{\bf Input:} a GLC $\kappa$ of \eqref{2_5}.
		\item{\bf Repeat} 
		\begin{enumerate}
			\renewcommand{\labelenumi}{\theenumi:}
			\item Measure the system state vector $x_k$.
			\item {\bf If} $k=0$,  then let $\mathbb{I}(x_k) = \mathbb{N}_{[1,n_c]}$, \\
			{\bf else} apply Algorithm \ref{alg:brief} with $\left(z^*(x_{k-1}), \mathbb{A}(x_{k-1})\right)$ to compute $\mathbb{I}(x_k)$. 
			\item Solve the trimmed QP($x_k$) with $\mathbb{I}(x_k)$ and record the results as $\left(z^*(x_k), \mathbb{A}(x_k)\right)$. \label{2_5_3}
			\item Apply the first input vector of $z^*(x_k)$ to the LTI system as \eqref{2_4}.
			\item $k\leftarrow k+1$.
		\end{enumerate} 
	\end{itemize}
\end{algorithm}

Now we are ready to quantify the number of linear inequalities $|\mathbb{I}(x_k)|$ in the trimmed QP($x_k$) to evaluate the performance of Algorithm \ref{alg:adaptive}.

\begin{theorem} \label{thm:adaptive}
	Suppose that Assumptions \ref{ass:LICQ} and \ref{ass:basic} hold. Consider Algorithm \ref{alg:adaptive}, the following statements are in force
	\begin{enumerate}[(a)]  \label{5_4}
		\item \label{5_4_2} If $k \ge K_i$, then $|\mathbb{I}(x_k)| \le n_z + i $.
		\item \label{5_4_3} If the origin is an interior of $\mathcal{X},\mathcal{U}$ and  $\mathcal{X}_N$, then $|\mathbb{I}(x_k)| = 0$ for all $k \ge \widehat{K}$, where $K_i$ and $\widehat{K}$ are given as
	\end{enumerate}
			\begin{eqnarray*}
				&&\hspace{-0.8cm} K_i = \left\lceil \log_\beta \frac{\sigma_{i}}{c\Vert x_0 \Vert (1 + \beta^{-1}) \sqrt{1 + \kappa^2}}\right\rceil, \ i \in  \mathbb{N}_{[1, n_c]} \\
				&& \hspace{-0.8cm}\widehat{K} = \max \{\widehat{K}_1, \widehat{K}_2\}, ~\text{and}~\label{5_5_2}  \\
				&&\hspace{-0.8cm} \widehat{K}_1 = \max_{j \in \mathbb{N}_{[1, n_c]}} \left\lceil \log_\beta \frac{w_j}{c \Vert x_0 \Vert \beta^{-1} (\kappa \Vert G_j \Vert + \Vert S_j \Vert)} \right\rceil, \\
				&&\hspace{-0.8cm} \widehat{K}_2 = \max_{j \in \mathbb{N}_{[1, n_c]}} {
					\left\lceil \log_\beta \frac{w_j}{\rho_j} \right\rceil},\\
						&&\hspace{-0.8cm}\rho_j={c \Vert x_0 \Vert \left( \kappa \Vert G_j \Vert (1 + \beta^{-1}) + \Vert S_j \Vert + \Vert  G_j H^{-1}F^T\Vert \beta^{-1}\right)}.
				\end{eqnarray*}	
			
	\end{theorem}
\begin{pf}
	See Appendix \ref{app:online_mpc}. 
\end{pf}

In view of Lemma \ref{lem:D}, the sequence $\{K_i\}$ is decreasing and thus the upper bound of the number of inequality constraints in QP($x_k$) decreases along with timestep. 
	
For stabilization problem with MPC,  it follows from \eqref{2_2} that $z^*(0)=0$. Jointly with \eqref{2_5_2}, it implies that $w > 0$ and thus $\widehat{K}$ is finite. Hence, the number of inequality constraints in the trimmed QP($x_k$) decreases to zero after a finite timestep $\widehat{K}$. 

\begin{remark}[Connection to the dual-mode MPC \cite{michalska2002robust}]
In the dual-mode MPC, a terminal set is precomputed offline such that within this set the control law \eqref{2_4} can be solved via the \emph{constraint-free} QP. By contrast, Algorithm 3 achieves this objective adaptively within the online closed-loop operation. Moreover, our approach does not require the computation of \emph{any} terminal set.
\end{remark}
\begin{remark}
	We emphasize that Algorithm \ref{alg:adaptive} does \emph{not} rely on Assumptions \ref{ass:LICQ} or \ref{ass:basic} to preserve the optimal solution and remove redundant constraints. The method remains valid in a wide range of application scenarios where these assumptions may fail, including trajectory tracking MPC and cases where the origin is not an interior point of all constraint sets. In such general settings, our trimmed MPC continues to preserve the optimal solution, though it may not eventually become completely constraint-free.
\end{remark}

\subsection{Learn to trim MPC from offline solved mp-QP instances} \label{sec:mpc_offline}
Of particular note is that our trimming algorithm cannot only reduce redundant inequalities from the previously solved instances in the online closed-loop system but also those in the offline design. To elaborate it, suppose that we have already solved $q$ instances of MPC {\em offline} and obtained that   
\begin{equation} \label{Ddef}
	\mathcal{D}= \{(\widehat{x}^i, z^*(\widehat{x}^i), \mathbb{A}(\widehat{x}^i))\}_{i \in \mathbb{N}_{[1,q]}}
\end{equation}
where $\widehat{x}^i$ corresponds to the parameter vector of an offline solved instance.

In the online closed-loop system, we select a triple $(\widehat{x}, z^*(\widehat{x}), \mathbb{A}(\widehat{x}))$ from \eqref{Ddef} with the parameter vector $\widehat{x}$ that is closest to $x_k$, i.e.,  
			 \begin{equation} \label{mindist}
			 \widehat{x}=\text{argmin}_{i\in\mathbb{N}_{[1,q]}}\Vert \widehat{x}^i- x_k\Vert
			 \end{equation}
instead of using all of them in Section \ref{sec:online_multiple}, and then apply Algorithm \ref{alg:brief} to trim mp-QP instances of MPC. 

\begin{remark}[Selection of parameter vectors of offline mp-QP instances] \label{rem_partition}By Theorem \ref{thm:distance}, the distance between parameter vectors of two instances is essential to reduce redundant inequalities. In this view, we should maximize the chance to find a nearby parameter vector of offline solved instance, which motivates to adopt the Voronoi diagram \cite{aurenhammer1991voronoi} to select parameter vectors. Specifically, they are chosen as the centroids of $q$-Voronoi cells $\{\mathcal{C}_i\}$ of the invariant set $\mathcal{X}$, i.e.,  
 \begin{equation}
 \label{Voronoi}
	\mathcal{C}_i = \{x \in \mathcal{X} \ | \ \Vert x - \widehat{x}^i \Vert \le \Vert x - \widehat{x}^j \Vert, \forall j \neq i \}, \ i \in \mathbb{N}_{[1, q]}.
 \end{equation}
That is,  Voronoi diagrams partition $\mathcal{X}$ into $q$-Voronoi cells, i.e., $\mathcal{X} = \cup_{i \in \mathbb{N}_{[1, q]}} \mathcal{C}_i$.  In this case, the time complexity of \eqref{mindist} is $\mathcal{O}(q)$. 

An alternative is to use a grid partition of $\mathcal{X}$, and let $\widehat{x}^i$ be the grid points of the partition. Then, the time complexity in \eqref{mindist} is $\mathcal{O}(\log(q))$.\qed
\end{remark}

\begin{figure}[t]
	\centering
        \includegraphics[width=8.5cm]{{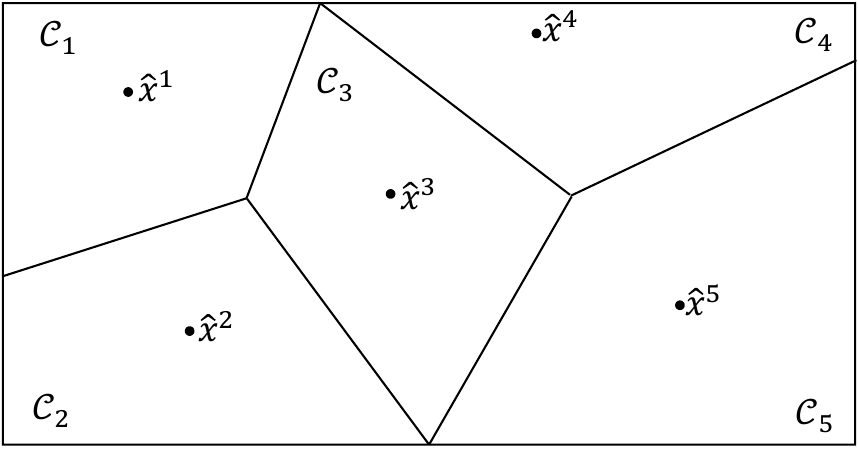}}
	\caption{Voronoi diagrams. The rectangular region represents the feasible set $\mathcal{X}$, which is partitioned into five Voronoi cells using \eqref{Voronoi}. }
	\label{fig:region}
\end{figure}

By Theorem \ref{thm:distance}, we can easily show the performance of our trimming method using the data of offline solved mp-QP instances in \eqref{Ddef}. Particularly, we define the maximum distance of  $\mathcal{X}$ as follows:
\begin{equation} \label{maxd}
	d(\mathcal{X},\mathcal{D}) = \max_{i\in\mathbb{N}_{[1,q]}} \max_{x \in \mathcal{C}_i} \Vert x - \widehat{x}^i \Vert.
\end{equation}
\begin{coro} \label{thm:mpc_offline}
	  Suppose that Assumption \ref{ass:LICQ} holds. If $d(\mathcal{X},\mathcal{D}) \le \sigma_i$ for some $i \in \mathbb{N}_{[1,n_c]}$,  then $\mathbb{I}(x_k)$ of the trimmed version of \eqref{2_2} with offline solved mp-QP instances in \eqref{Ddef} satisfies that $$\mathbb{I}(x_k) \le n_z + i~\text{for any}~ k \in \mathbb{N}.$$ 
\end{coro}

Clearly, the more offline solved instances used, the smaller the distance $d(\mathcal{X},\mathcal{D})$, which potentially leads to a slimmer mp-QP for the trimmed version. 

\begin{remark}[Hybrid learning]\label{rem_hybrid}
	Our trimming method can also easily take advantages of both offline and online solved instances, e.g., one can jointly use the previously solved instance at one timestep lag of the closed-loop system and the closest instance in $\mathcal{D}$. See Section \ref{sec:online_multiple} for details of learning from multiple solved instances. 
\end{remark}

\begin{remark}[Comparison with the explicit MPC] The control law of the explicit MPC \cite{alessio2009survey} takes a piecewise affine structure in the parameter vector and can be obtained via solving mp-QP offline. To implement it online, one simply adopts the table-lookup method. However, the number of table cells grows exponentially with the size of mp-QP in the worst case, and thus is only well suited for mp-QP with a moderate size. For example, we cannot implement the explicit MPC on our personal laptops for the simulation in Section \ref{sec:sim}. In comparison, our trimmed MPC is though not fully solved offline, it easily takes benefits from offline solved instances to remove redundant constraints of instances in the closed-loop system, based on which the trimmed instances can be much cheaper to solve online.  Importantly, it is quite flexible and easy to take advantages of both offline and online solved instances. 
\end{remark}

\section{Simulation} \label{sec:sim}
In this section, we evaluate the performance of the proposed method through simulations on two distinct examples. The first is an oscillating mass system, which serves as a typical example for stabilizing control. The second is an unmanned aerial vehicle (UAV) trajectory tracking task, representing a highly dynamic scenario where the active set changes frequently. All simulations were implemented in Python and executed on a standard laptop equipped with an Intel Core i5-11400F processor and 16 GB of RAM.

\begin{figure}[t]
	\centering
	\includegraphics[width=8.5cm]{{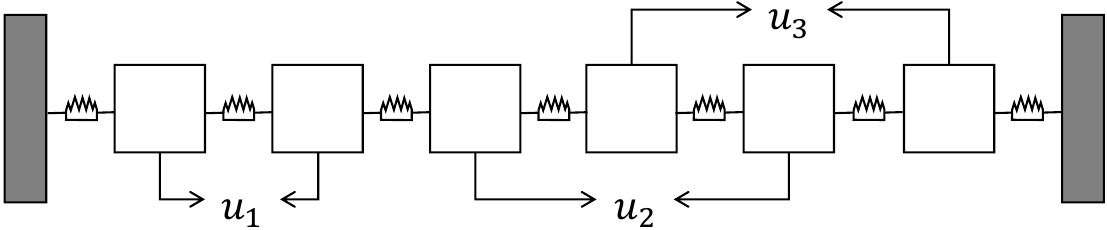}}
	\caption{An oscillating mass system. The arrow denotes the direction of the force that is exerted to masses.}
	\label{fig:mass}
\end{figure}

\subsection{Application to Oscillating Mass System}

The oscillating mass system \cite[Section V-A]{wang2009fast} consists of six masses with unit value that are connected by springs (with the spring constants of $1$). See Fig \ref{fig:mass} for an illustration. Let the state vector be the position and velocity of all masses, i.e., $x \in \mathcal{R}^{12}$.   Three actuators exert tensions to masses under  a maximum control force of $\pm 0.5$, and the displacements of masses cannot exceed $\pm 4$. 

For the MPC design of \eqref{2_1}, we use the ZOH method to discretize the continuous-time system at a sampling rate $10Hz$, and adopt the procedures in \cite[Chapter 2]{rawlings2017model} to satisfy Assumption \ref{ass:basic}. In particular, let the weighting matrices $Q, R$ be simply set as identity matrices, and $P$ in the terminal cost is set as the positive definite solution to the discrete-time algebraic Riccati equation, i.e., $P = A^TPA - A^TPB(R + B^TPB)^{-1} B^TPA + Q$. Let the associated optimal gain be $K^*= -(R + B^TPB)^{-1}B^TPA$, based on which we set $\mathcal{X}_N$ to be the maximal invariant set of the closed-loop system $x_{k+1} = (A + BK^*)x_k$ subject to $\|x_k\|_\infty\le 4$ and $\|u_k\|_\infty=\|K^*x_k\|_\infty\le 0.5$. Specifically, we use \cite{gilbert1991linear} to compute $\mathcal{X}_N$, which results in $450$ linear inequalities. One can easily verify that the mp-QP instance in \eqref{2_2} is recursively feasible starting from any $x_0 \in \mathcal{X}_N$. Moreover, we randomly select $20$ initial state vectors from $\mathcal{X}_N$, each of which is run for $100$ timesteps in the closed-loop systems under the control law in \eqref{2_4}, and report averaged results. We set the prediction horizon to $N = 50$, resulting in mp-QP dimensions of $n_x = 12$, $n_z = 50\times 3=150$, and $n_c = 50 \times (12 + 6) + 450 = 1350$.

\textbf{Comparison with warm-started MPC and ca-MPC.}
We compare our trimmed MPC (Algorithm \ref{alg:adaptive}) against three baselines: the standard MPC, the ca-MPC proposed in Nouwens et al. \cite{nouwens2023constraint}, and the warm-started MPC. In the warm start, the working set from the previous solution is shifted in time to serve as the initial guess. For Algorithm \ref{alg:adaptive}, we utilize NumPy library for vectorized matrix operations in the constraint removal, and precompute the GLC constant $\kappa = 67.24$ via Lemma \ref{lem:LipschitzPro} using the commercial solver {\it Gurobi}. All the QP instances are solved using DAQP \cite{arnstrom2022dual}, a state-of-the-art active-set solver.

\begin{figure}[t]
	\centering
        \includegraphics[width=8.5cm]{{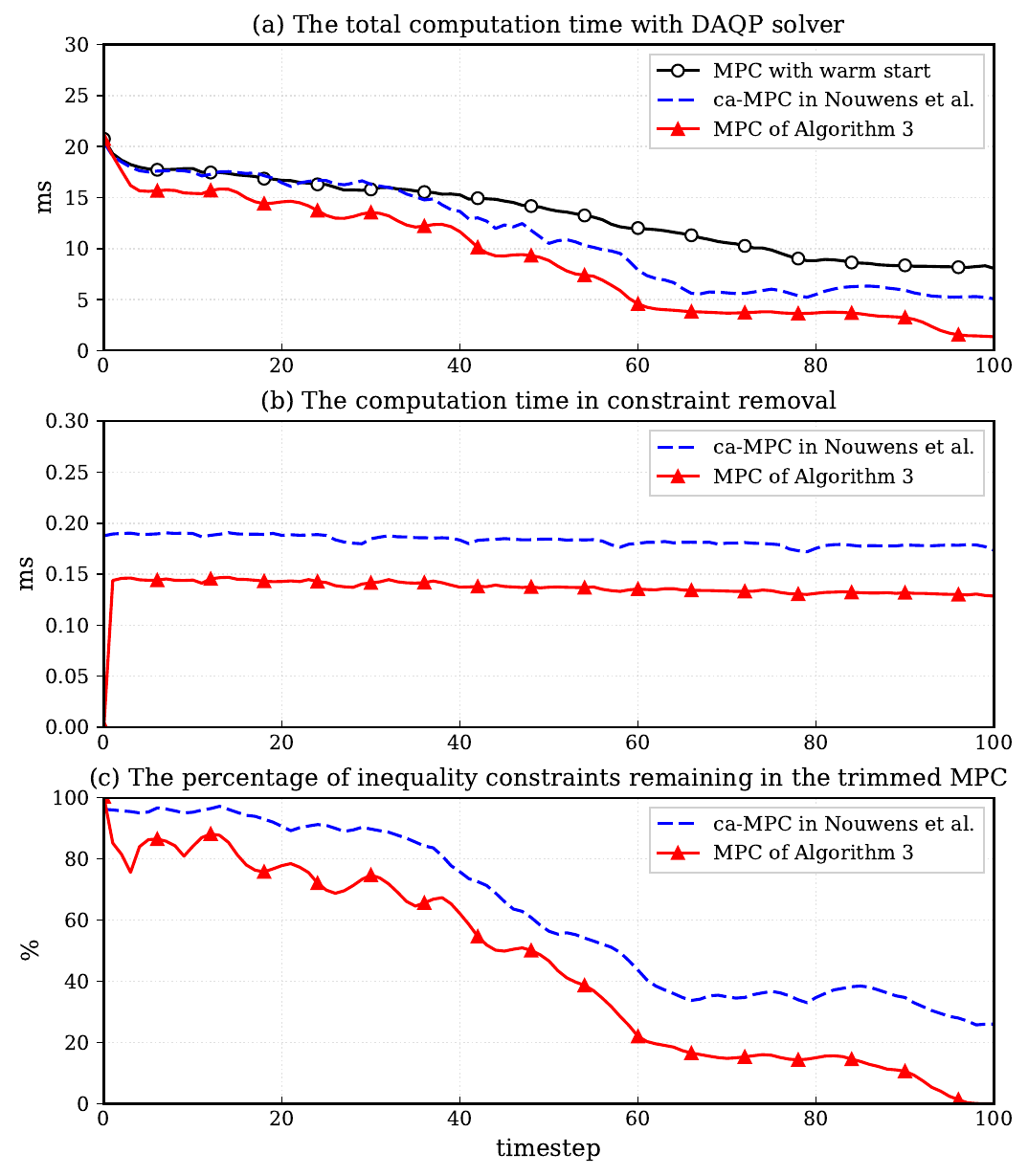}}
	\caption{Performance comparison between MPC of the proposed Algorithm 3, warm-started MPC, and MPC with constraint removal in \cite{nouwens2023constraint}. The top subfigure illustrates the total computation time, including both the constraint removal and the QP solution times.}
	\label{fig:adaptive}
\end{figure}

By setting the standard MPC as a baseline, we record the averaged percentages of computation time and constraints in Fig. \ref{fig:adaptive}, which indeed confirms superiorities of our trimming MPC. The observed performance gain originates from our approach's ability to simplify the mp-QP instance independently of the solver by removing provably redundant constraints beforehand. This reduction is particularly significant for active-set solvers, whose worst-case complexity grows exponentially with the number of constraints \cite{wright2006numerical}. Thus, our method serves as a complementary preprocessing step, and eventually leading to the constraint-free results observed in Algorithm \ref{alg:adaptive}.

\begin{figure}[t]
	\centering
	\includegraphics[width=8.5cm]{{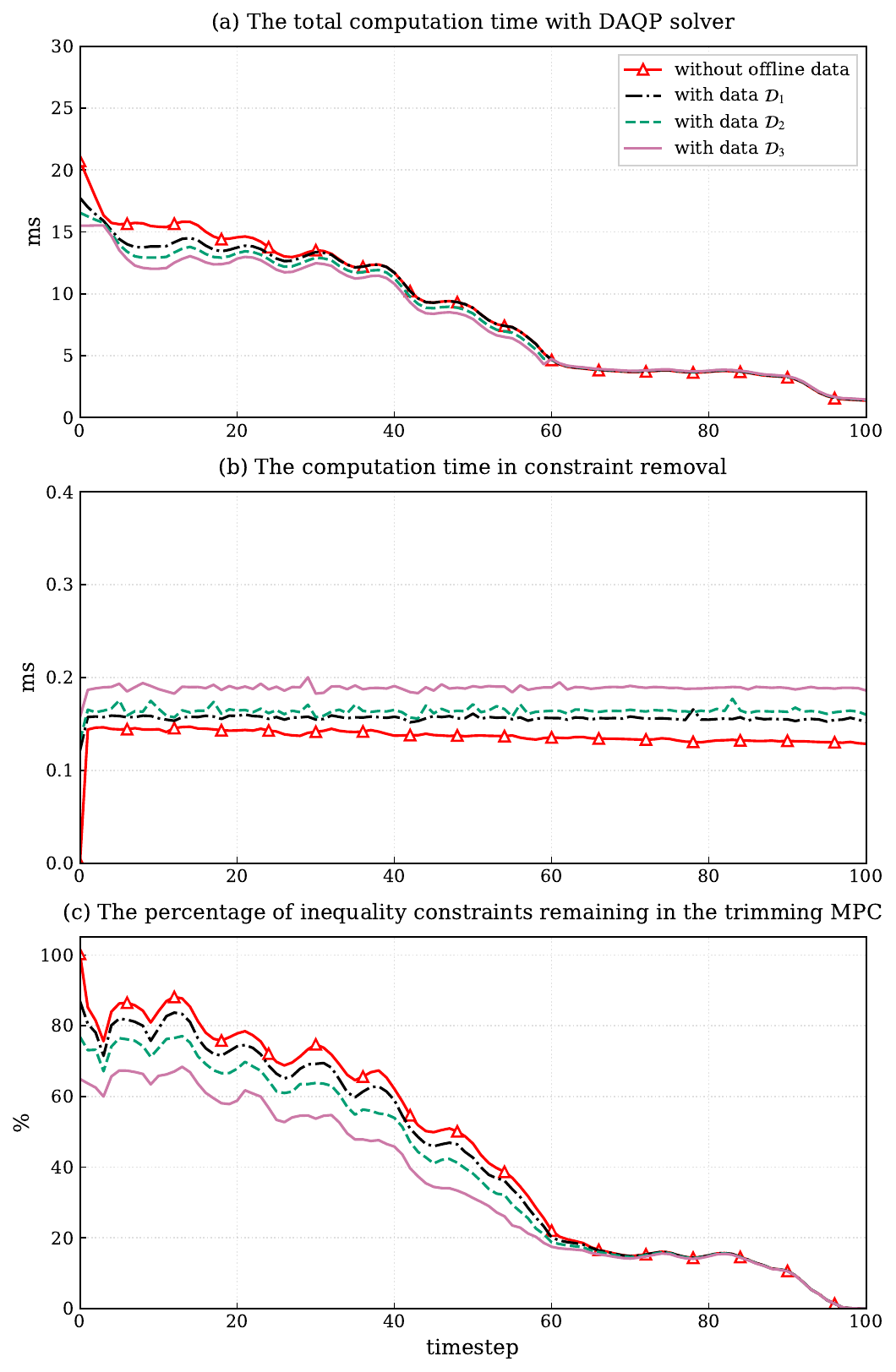}}
	\caption{Performance comparison between warm-started MPC and the trimmed MPC with different offline/online data.}
	\label{fig:per_constraint}
\end{figure}

\textbf{Comparison with the explicit MPC.} We use the MPT3 Toolbox \cite{MPT3} to compute the explicit MPC. Unfortunately, the algorithm has run for more than $72$ hours, finds more than $51954$ regions and still does not terminate.  That is, we are essentially unable to use the explicit MPC for this case in our laptop.  To illustrate our results in Section \ref{sec:mpc_offline},  we adopt the grid partition method in Remark \ref{rem_partition} over the invariant set $\mathcal{X}_N$ and construct three data sets $\mathcal{D}_1, \mathcal{D}_2~\text{and}~\mathcal{D}_3$, whose details are reported in Table \ref{tab:D}. Since the LICQ condition in  Assumption \ref{ass:LICQ} holds, we jointly learn from both offline and online solved instances as stated in Remark  \ref{rem_hybrid}. As shown in Fig. \ref{fig:per_constraint}, the effect of offline solved instances decreases in timestep, and the online solved instances dominate the importance of removing constraints. 
\begin{table}[t] 
	\centering 
	\caption{ Information of the data set}
	\label{tab:D}
	\begin{tabular}{c l c l c}\hline 
		Data set &   & Coordinate distance &   & Number of data \\
		\hline
		$\mathcal{D}_1$ & & 0.3 &  & 86 \\ 
		$\mathcal{D}_2$ & & 0.25 &  & 312 \\ 
		$\mathcal{D}_3$ & & 0.2 &  & 1195 \\ 
		\hline
	\end{tabular}	 
\end{table}

\begin{figure}[t]
	\centering
	\includegraphics[width=8.5cm]{{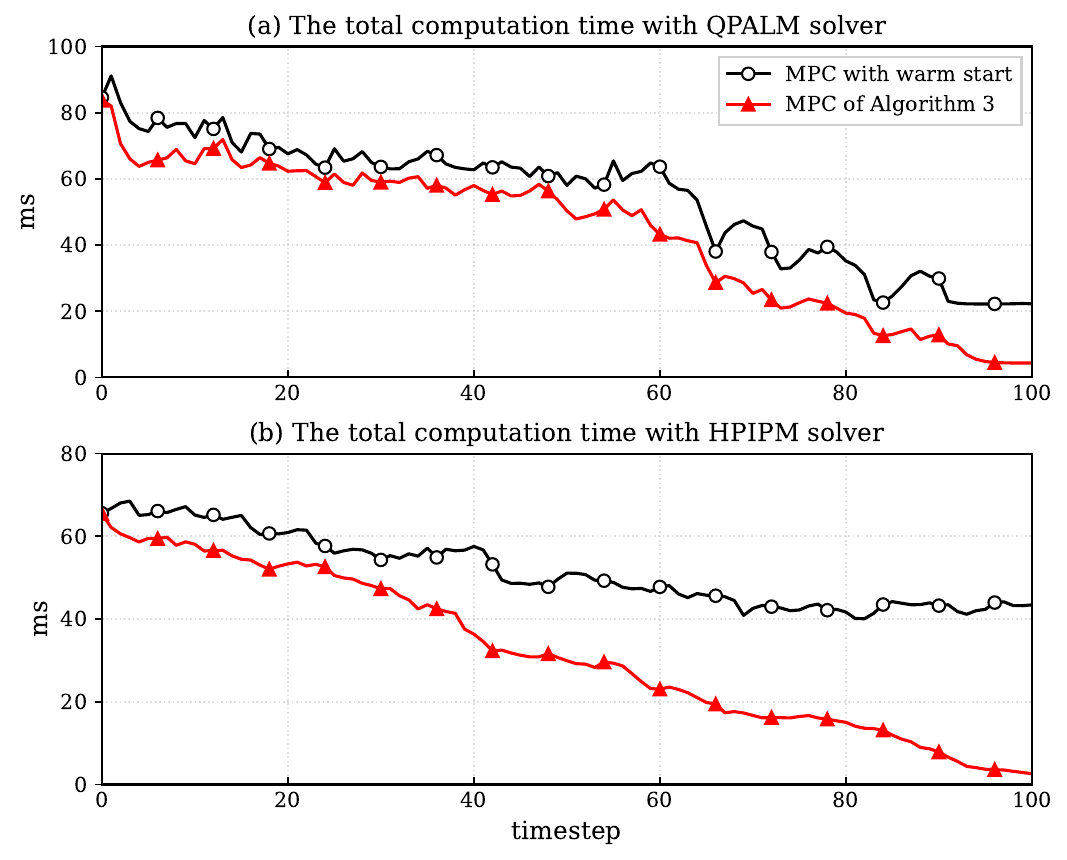}}
	\caption{Comparison of the average computation time under two QP solvers.}
	\label{fig:solver}
\end{figure}

\textbf{Comparison of the computation time under different solvers.} We further compare the trimmed MPC in Algorithm \ref{alg:adaptive} with the warm-started MPC under the different solvers. This simulation employs two distinct QP solvers: the interior-point method based HPIPM \cite{FRISON20206563} and the augmented Lagrangian based QPALM \cite{hermans2019qpalm}. Specifically, the warm start strategy is implemented by shifting the optimal primal and dual solutions from the previous time step to initialize the current optimization. The average computation time required by each method under both solvers at each time step is depicted in Fig \ref{fig:solver}. The result demonstrates that the acceleration achieved by the proposed method is not specific to a single solver but is indeed effective across different QP solvers.

\subsection{Application to UAV Trajectory Tracking}
For the trajectory tracking scenario, we apply Algorithm 3 to control a UAV, where the model and the MPC formulation are adopted from \cite{kamel2017linear}. Specifically, the model is linearized around the hovering condition. We specially adopt stricter bound settings for the roll and pitch angles, where the bounds are restricted to approximately 0.38 rad to encourage the active set to change more dynamically during tracking. Although the reference trajectory is time-varying, the problem can still be formulated as a mp-QP by simply including the reference trajectory in the parameter vector $x$. With the prediction horizon set to 20, this formulation results in mp-QP dimensions of $n_x = 8$, $n_z = 20 \times 3 = 60$, and $n_c = 20 \times 6 = 120$. The GLC constant is computed as $\kappa = 32.11$, and all the instances are solved using DAQP. The simulation is conducted for 400 timesteps at a frequency of 10 Hz.

\begin{figure}[t]
	\centering
	\includegraphics[width=8.5cm]{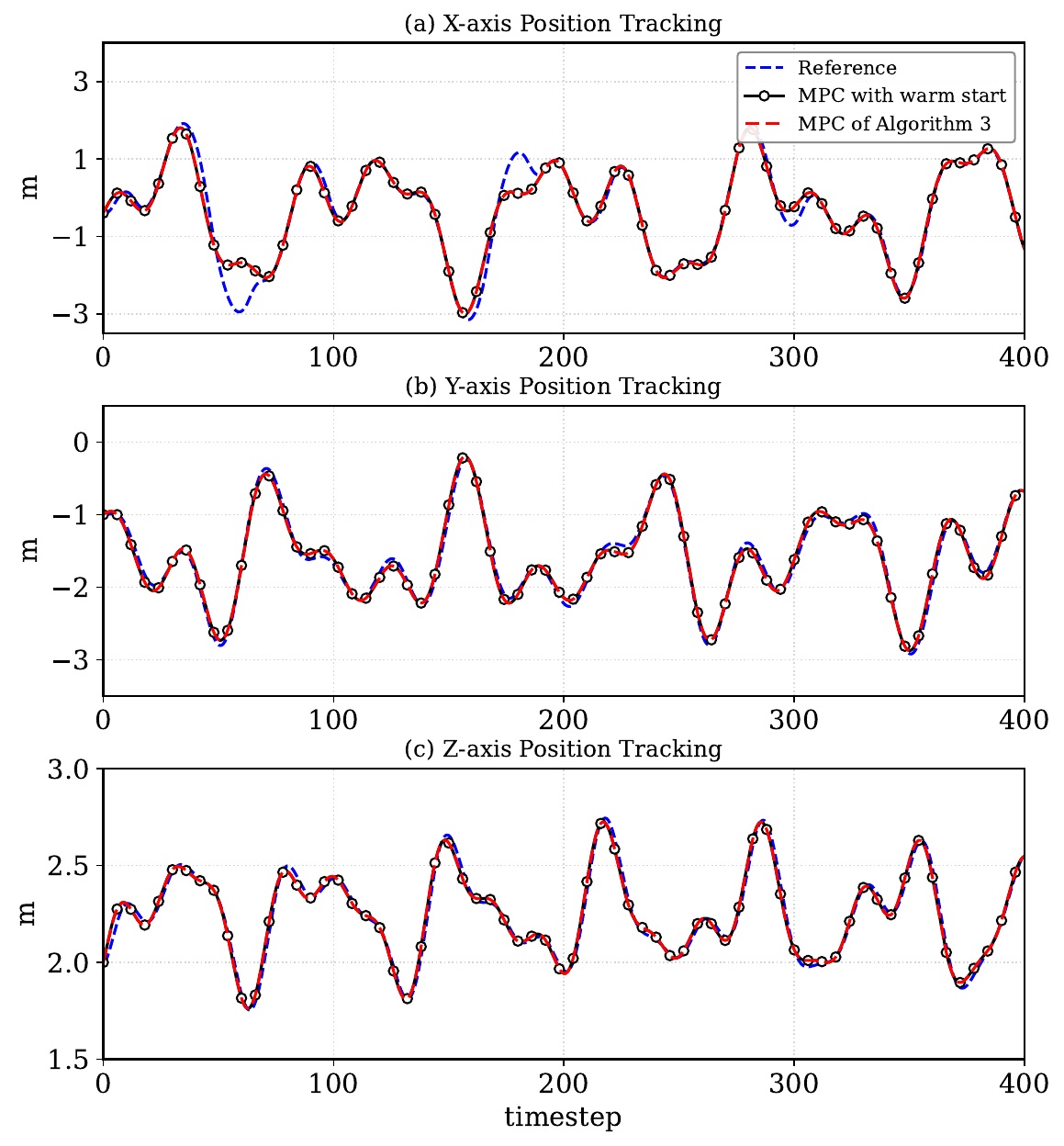}
	\caption{3D trajectory tracking performance of warm-started MPC and the MPC of Algorithm 3 for the UAV.}
	\label{fig:uav_tracking}
\end{figure}

\begin{figure}[t]
	\centering
	\includegraphics[width=8.5cm]{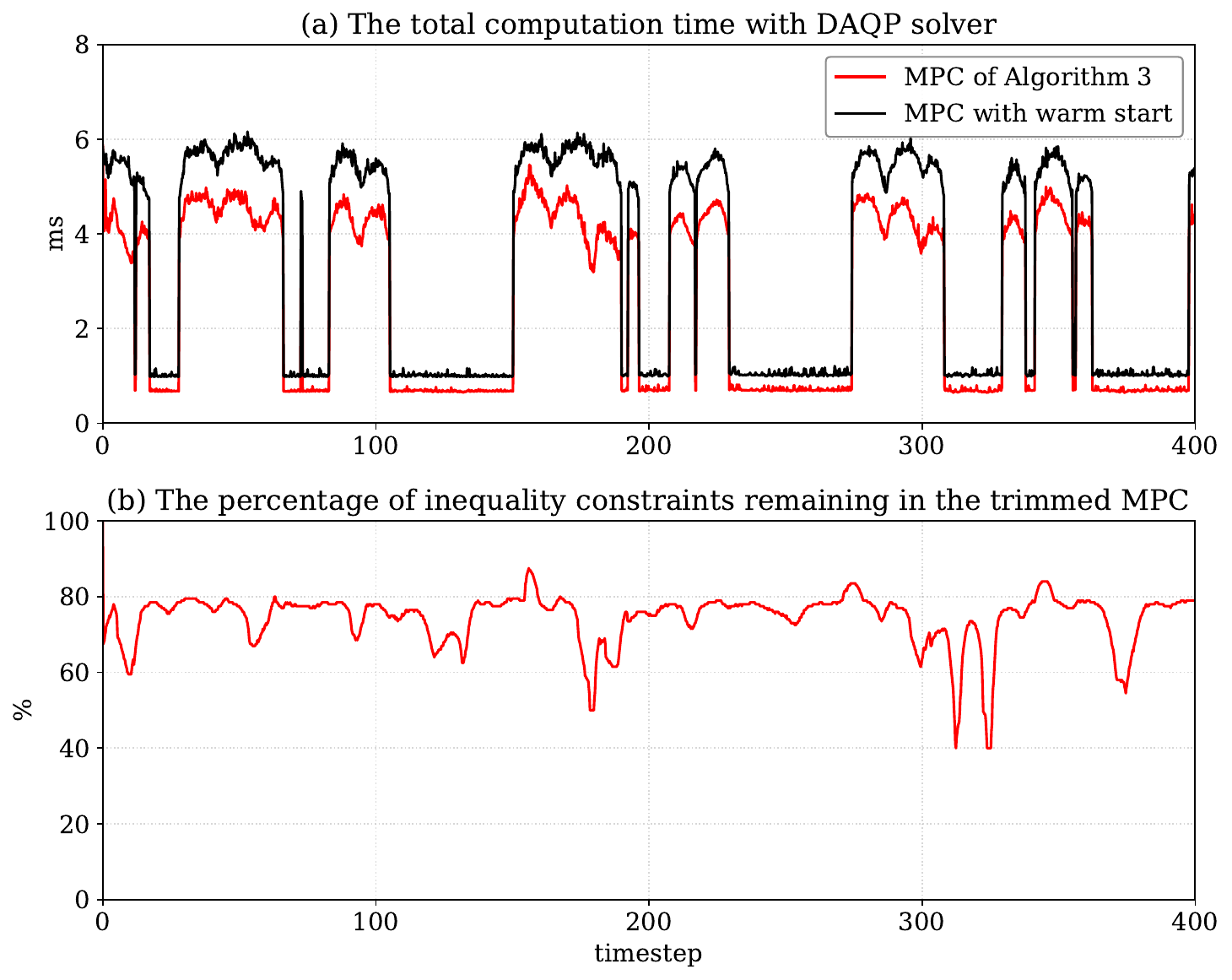}
	\caption{Performance comparison between warm-started MPC and the MPC of Algorithm 3 during the UAV trajectory tracking.}
	\label{fig:uav_online}
\end{figure}

The results are presented in Fig. \ref{fig:uav_tracking} and Fig. \ref{fig:uav_online}. Fig. \ref{fig:uav_tracking} illustrates the 3D trajectory tracking curves of the UAV. As shown, the proposed method achieves a tracking performance that is completely identical to warm-started MPC. This confirms that our proposed method strictly preserves the optimal solution without any performance degradation, even in a trajectory tracking environment. Moreover, as depicted in Fig. \ref{fig:uav_online}, our method successfully trims approximately 30\% of constraints at each step, which directly leads to a noticeably reduced computation time. We can also observe that the computation time occasionally drops to a very small value. This simply occurs because the active set changes frequently during tracking; when the active set becomes empty, the solver can find the optimal solution extremely fast.

\section{Conclusion} \label{sec:conclusion}
In this work, we have presented novel methods to learn from the results of previously solved mp-QP instances to accelerate the unsolved instances, which are reliably trimmed via constraint removal. Then, we extended to trim mp-QP instances of MPC by utilizing both online and offline solved instances. Numerical results were included to support our theoretical findings. Overall, learning to remove redundant constraints appears to be a promising technique to solve a sequence of multiparametric programs, which will be further exploited to accelerate linear time-varying MPC in our future works. Additionally, scaling the offline computation of the global Lipschitz constant to high-dimensional MPC applications represents another important avenue for future research.

\appendix 
For writing brevity, we introduce some notations in this section. $\mathbb{B}=\{0,1\}$. For any matrix $A \in \mathcal{R}^{m \times n}$ and any set $\mathbb{I} \subseteq \mathbb{N}_{[1,m]}$, let $A_{\mathbb{I}}$ denote the submatrix that corresponds to the rows indexed by $\mathbb{I}$. For two vectors $a$ and $b$ of the same dimension, we denote their Hadamard product, i.e., pairwise product, by $a \circ b$. For matrices $A$ and $B$, we denote their Kronecker product by $A \otimes B$. The operator $\text{col}(·)$ stacks its arguments in column vectors or matrices of compatible dimensions.
\subsection{Proof of Lemma \ref{lem:Lipschitz}} \label{app:Lipschitz}

The KKT conditions of the trimmed mp-QP of \eqref{2_6} with $\mathbb{I}(x) = \mathbb{I}$ \cite[Section 3.4]{bertsekas2016nonlinear}  are given as
\begin{equation} 
	\begin{aligned}
		Hz + F^Tx + G_\mathbb{I}^T \lambda = 0, \\
		\lambda \circ (G_\mathbb{I}z - w_\mathbb{I} - S_\mathbb{I} x) = 0, \\
		\lambda \ge 0, \\
		G_\mathbb{I}z \le w_\mathbb{I} + S_\mathbb{I}x.
	\end{aligned}  \label{A_0}
\end{equation}
where  we denote the primal solution of \eqref{A_0} as $z^*(x, \mathbb{I})$.

Consider the following equations:
\begin{subequations} \label{A_1}
	\begin{align} 
		Hz + F^Tx + G^T \lambda = 0, \label{A_1_1} \\
		\lambda \circ (Gz - w - S x) = 0, \label{A_1_2} \\
		\lambda \ge 0. \label{A_1_3} 
	\end{align}
\end{subequations}

First, we show that the set of primal solutions of \eqref{A_1} can be given as 
\begin{equation}\{z^*(x, \mathbb{I}) \ | \forall \mathbb{I} \subseteq \mathbb{N}_{[1,n_c]} \}\label{soluset}
\end{equation}
On the one hand, for any $z^*(x,\mathbb{I})$ of \eqref{soluset}, there exists a dual vector $\lambda^*(x, \mathbb{I})$ such that  the pair of primal-dual vectors $(z^*(x,\mathbb{I}), \lambda^*(x, \mathbb{I}))$ solves  \eqref{A_0}.  Then, we augment $\lambda^*(x, \mathbb{I})$ by padding zero elements to those positions in the complement of  $\mathbb{I}$, and obtain a new vector $\lambda^* \in \mathcal{R}^{n_c}$, i.e., $\lambda_{\mathbb{I}}^*= \lambda^*(x, \mathbb{I})$ and $\lambda_{\mathbb{I}^c}^* = 0$. One can validate that the pair of primal-dual vectors $(z^*(x, \mathbb{I}),\lambda^*)$ satisfies \eqref{A_1}.  On the other hand, let $(z^*, \lambda^*)$ be a pair of primal-dual vectors of \eqref{A_1}. Define the index set $\mathbb{I} = \{i \in \mathbb{R}_{[1,n_c]}\ | \ G_iz^* \le w_i + S_i x \}$. That is, $(z^*, \lambda_\mathbb{I}^*)$ satisfies \eqref{A_0}, and $z^* = z^*(x, \mathbb{I})$ is an element of  \eqref{soluset}. 

The rest of proof is similar to \cite[Section V]{9928332}. We rewrite \eqref{A_1} as follows:
\begin{subequations} \label{A_2}
	\begin{align} 
		Hz + F^Tx + G^T \lambda &= 0, \label{A_2_1} \\
		Gz - w - S x + r &= 0, \label{A_2_2} \\
		-\overline{r} (1 - \delta) \le r &\le \overline{r} (1 - \delta), \label{A_2_3} \\
		0 \le \lambda &\le \overline{\lambda} \delta. \label{A_2_4} 
	\end{align}
\end{subequations}
where slack vector $r \in \mathcal{R}^{n_c}$, binary vector $\delta \in \mathbb{B}^p$, $\overline{r}$ and $\overline{\lambda}$ are ``big'' values, i.e., upper bounds of $r$ and $\lambda$, respectively.

Under the fixed $\delta$, deriving the expression of $z$ about $x$ from \eqref{A_2} yields affine function $z(x;\delta) = K(\delta)x + c(\delta)$. To express $K(\delta)$, we perturb the parameter vector $x$ along the canonical basis vectors in $\mathcal{R}^n$. Specifically, we introduce the auxiliary variables $\{x_i, z_i, r_i, \lambda_i\}$ for $i = 1,2,\ldots,n_x$ and additional mixed integer linear constraints as 
\begin{subequations} \label{A_3}
	\begin{align} 
		x_i = x + e_i, \label{A_3_0} \\
		Hz_i + F^Tx_i + G^T \lambda &= 0, \label{A_3_1} \\
		Gz_i - w - S x_i + r_i &= 0, \label{A_3_2} \\
		-\overline{r} (1 - \delta) \le r_i &\le \overline{r} (1 - \delta), \label{A_3_3} \\
		-\overline{\lambda_i} \delta \le \lambda_i &\le \overline{\lambda_i} \delta. \label{A_3_4} 
	\end{align}
\end{subequations}
where $e_i \in \mathcal{R}^n$ is the $i$-th vector of the canonical basis.

It holds that $z_i = K(\delta)x_i + c(\delta)$, and hence that:
\begin{equation} \notag 
	\begin{aligned}
		[z_1, z_2, \ldots, z_{n_x}]  &= K(\delta)[x_1, x_2,\cdots,x_{n_x}] + c(\delta) \\
		& = K(\delta) (x \otimes \textbf{1}^T + I) + c(\delta) \\
		& = z \otimes \textbf{1}^T + K(\delta)
	\end{aligned}
\end{equation}
Then we turn the linear segment $K(\delta)$ into the left side:
\begin{equation} \label{A_4}
	K(\delta) = [z_1, z_2, \ldots, z_{n_x}] - z \otimes \textbf{1}^T 
\end{equation}

We next compute the $l_1$-norm $\|K(\delta)\|_1$. To ensure completeness, we reproduce the derivation based on \cite{9928332}. The $l_1$-norm computation can be formulated as the following linear programming problem:
\begin{equation}
	\|K(\delta)\|_1 = \min_{\eta} \ c^\top \eta \quad \text{s.t.} \quad M\eta \le h(\delta), \label{A_6}
\end{equation}
where $\eta \in \mathbb{R}^{mn+1}$ is the decision variable, $c = \text{col}(\mathbf{0}, 1)$, and the constraint parameters are defined as:
\begin{equation}
	M := \begin{bmatrix}
		\mathbf{1}^\top \otimes I & -\mathbf{1} \\
		-I & \mathbf{0} \\
		-I & \mathbf{0}
	\end{bmatrix}, \quad h(\delta) := \begin{bmatrix}
		\mathbf{0} \\
		\text{vec}(K(\delta)) \\
		-\text{vec}(K(\delta))
	\end{bmatrix}.
\end{equation}

Since the optimization problem \eqref{A_6} is always feasible, strong duality holds. Now consider its associated dual problem:
\begin{equation}
	\|K(\delta)\|_1 = \max_{\mu \ge 0} \ -h(\delta)^\top \mu \quad \text{s.t.} \quad M^\top \mu = -c. \label{A_8}
\end{equation}
As established in \cite[Lemma 5.1]{9928332}, this dual problem admits a binary optimizer. This property allows us to restrict the search space to binary variables without loss of optimality. Thus, we replace the constraint $\mu \ge 0$ with $\mu \in \mathbb{B}^{n_\mu}$, where $n_\mu = n(1+2m)$.

The objective function $-h(\delta)^\top \mu$ in \eqref{A_8} contains bilinear terms. To linearize this, let $k(\delta) := \text{vec}(K(\delta))$ and partition the binary vector as $\mu = \text{col}(\mu^0, \mu^1, \mu^2)$ with $\mu^0 \in \{0,1\}^n$ and $\mu^1, \mu^2 \in \{0,1\}^{mn}$. The objective function simplifies to:
$$
-h(\delta)^\top \mu = -k(\delta)^\top \mu^1 + k(\delta)^\top \mu^2.
$$
We introduce auxiliary variables $y \in \mathbb{R}^{2mn}$ to represent the linearized terms corresponding to $k(\delta)^\top \mu^1$ and $-k(\delta)^\top \mu^2$. The logical implications (e.g., $\mu^1_j=1 \Rightarrow y^1_j = k(\delta)_j$ and $\mu^1_j=0 \Rightarrow y^1_j = 0$) are enforced via Big-M constraints:
\begin{equation} \label{A_9}
	E \ \text{col}(k(\delta), \mu, y) \le g.
\end{equation}
where the matrix $E \in \mathbb{R}^{8mn \times (mn + n_\mu + 2mn)}$ and vector $g \in \mathbb{R}^{8mn}$ are constructed as:
\begin{equation}  \label{A_10}
	E := \begin{bmatrix}
		\mathbf{0} & -tX& I_{2mn} \\
		\mathbf{0} & -tX & -I_{2mn} \\
		-Y & tX & I_{2mn} \\
		Y & tX & -I_{2mn}
	\end{bmatrix}, \quad
	g := \begin{bmatrix}
		\mathbf{0} \\
		\mathbf{0} \\
		t\mathbf{1} \\
		t\mathbf{1}
	\end{bmatrix},
\end{equation}
where $X = [\mathbf{0}_{2mn \times n} \mid I_{2mn}]$, $Y = [I_{mn}; -I_{mn}]$. Here, the notations $[\cdot \mid \cdot]$ and $[\cdot; \cdot]$ denote horizontal and vertical matrix concatenations, respectively. $t$ is a sufficiently large constant \cite{9928332}.

Finally, combining the constraints \eqref{A_2}-\eqref{A_4} with the dual formulation \eqref{A_8}-\eqref{A_10}, we obtain the final MILP:
\begin{equation} \label{A_11}
	\begin{aligned}
		\Vert K(\delta) \Vert_1 = \max_{\delta, \mu, y} \quad & -\mathbf{1}^\top y \\
		\text{s.t.} \quad & \eqref{A_2}-\eqref{A_4}, \\
		& M^\top \mu = -c, \\
		& E \text{col}(\text{vec}(K(\delta)), \mu, y) \le g.
	\end{aligned}
\end{equation}
The optimal value of this MILP yields $\kappa_1$. The derivation for $\kappa_\infty$ follows an analogous procedure.

\subsection{Proof of Theorem \ref{pro:mul}} \label{app:mul}
The resultant index set of Algorithm \ref{alg:brief_mul} satisfies that:
\begin{equation}
	\mathbb{I}(x) = \cap_{k=1}^q \mathbb{I}^k(x), \label{app_B_1} 
\end{equation}
where $\mathbb{I}^k(x) = \mathbb{A}(\widehat{x}^k) \cup \mathbb{R}(x, \widehat{x}^k)$. 

It follows from KKT conditions of \eqref{2_5} and \eqref{2_6} with $\mathbb{I}(x) = \mathbb{I}^k(x)$ that
\begin{subequations}
	\begin{align}
		& Hz^*(x) + F^Tx + G^T \lambda^*(x) = 0, \notag \\
		& Hz^*(x, \mathbb{I}^k(x)) + F^Tx + G^T_{\mathbb{I}^k(x)} \lambda^*(x,\mathbb{I}^k(x)) = 0, \notag
	\end{align}
\end{subequations}
where $\lambda^*(x), \lambda^*(x,\mathbb{I}^k(x))$ are the Lagrange multipliers of the mp-QP \eqref{2_5} and the trimmed mp-QP \eqref{2_6} with the index set $\mathbb{I}^k(x)$, respectively. Since Theorem \ref{thm:data} ensures that $z^*(x, \mathbb{I}^k(x)) = z^*(x)$, it follows that
\begin{equation}
	G^T \lambda^*(x) = G^T_{\mathbb{I}^k(x)} \lambda^*(x,\mathbb{I}^k(x)). \label{app_B_2} 
\end{equation}

Clearly, $\lambda^*(x)$ is non-negative, and we denote the strict positive elements by $\underline{\mathbb{A}}(x) = \{j \in \mathbb{N}_{[1,n_c]} \ | \ \lambda_j^*(x) > 0 \}$. Let $\mathbb{A}(x, \mathbb{I}^k(x))$ denote the set of indices of active constraints in the trimmed mp-QP \eqref{2_6} with the index set $\mathbb{I}^k(x)$. Moreover, in view of the complementarity slackness \cite{bertsekas2016nonlinear}, $\lambda_j^*(x,\mathbb{I}^k(x)) = 0$ for any $j \notin \mathbb{A}(x, \mathbb{I}^k(x))$. It follows from \eqref{app_B_2} that
\begin{equation}
	\sum_{j \in \underline{\mathbb{A}}(x)} G^T_j \lambda_j^*(x) = \sum_{j \in \mathbb{A}(x, \mathbb{I}^k(x))} G^T_j  \lambda_j^*(x,\mathbb{I}^k(x)). \label{app_B_3} 
\end{equation}
Note that both $\underline{\mathbb{A}}(x)$ and $\mathbb{A}(x, \mathbb{I}^k(x))$ are subsets of $\mathbb{A}(x)$. Under Assumption \ref{ass:LICQ}, \eqref{app_B_3} implies that $\underline{\mathbb{A}}(x) \subseteq \mathbb{A}(x, \mathbb{I}^k(x))$. Then,  $\underline{\mathbb{A}}(x) \subseteq \mathbb{I}^k(x)$. It follows from \eqref{app_B_1} that $\underline{\mathbb{A}}(x) \subseteq \mathbb{I}(x)$. 

Since $\underline{\mathbb{A}}(x) \subseteq \mathbb{I}(x)$, $z^*(x)$ and $\lambda_{\mathbb{I}(x)}^*(x)$ satisfy the KKT conditions of the mp-QP \eqref{2_6}. Thus, $z^*(x)$ is also the optimal solution of the mp-QP \eqref{2_6}. By the uniqueness of its optimal solution, it follows that $
z^*(\widehat{x}, \mathbb{I}(x)) = z^*(\widehat{x})$ which completes the proof.

\subsection{Proof of Lemma \ref{lem:extended}} \label{app:extended}
(\ref{4_10_3}) Since $z^*(\widehat{x})$ is the optimal solution of mp-QP in \eqref{2_5}, $z^*(\widehat{x}) \in \mathcal{Z}(\widehat{x})$, i.e., $Gz^*(\widehat{x}) \le W + S \widehat{x}$. This yields $v(\widehat{x}) \in \mathcal{V}$.

(\ref{4_10_4}) Consider any $j \in \mathbb{E}(x)$.  For any $ q \in \mathcal{B}(z^*(\widehat{x}), \kappa \Vert x - \widehat{x} \Vert)$, let $v_1 = [x', q']'$. It holds that $ \Vert v_1 - v(\widehat{x}) \Vert \le \sqrt{1 + \kappa^2} \Vert x - \widehat{x} \Vert$. It follows from \eqref{4_10_2} that $H_j v_1 \le w_j$. Rewriting it yields that $Gq \le W + Sx$, i.e., $q \in \mathcal{Z}_j(x)$. This implies that $\mathcal{B}(z^*(\widehat{x}), \kappa \Vert x - \widehat{x} \Vert) \subseteq \mathcal{Z}_j(x)$. It follows from \eqref{3_8} that $j \in \mathbb{D}(x)$. Thus, $\mathbb{E}(x) \subseteq \mathbb{D}(x)$.

\subsection{Proof of Lemma \ref{lem:D}} \label{app:D}
The proof is organized into the following four parts.

(a) Existence of the minimum in \eqref{4_11_2}: 
Define the function
\begin{equation}
	S(v, r) = | \{ j \in \mathbb{N}_{[1, n_c]} \ | \ \mathcal{B}(v,r) \subseteq \mathcal{V}_j \}|.
\end{equation}
Clearly, for any $r \ge 0$, the set $\{S(v,r) \ | \ v \in \mathcal{V} \}$ is a subset of the finite integer set $\mathbb{N}_{[0, n_c]}$. Since every finite non-empty subset of integers attains its minimum, the minimum in \eqref{4_11_2} exists.

(b) Existence of the maximum in \eqref{4_11}: We introduce two step functions:
\begin{equation} \notag
	\epsilon(x) := \left\{ \begin{aligned} & 1, \quad x \ge 0 \\ & 0, \quad x < 0 \end{aligned} \right., \quad
	\delta(x) := \left\{ \begin{aligned} & 1, \quad x > 0 \\ & 0, \quad x \le 0 \end{aligned} \right..
\end{equation}
Note that $\epsilon(x) + \delta(-x) = 1$ holds for any $x \in \mathbb{R}$.
According to Remark \ref{rem_formula}, the condition $\mathcal{B}(v,r) \subseteq \mathcal{V}_j$ holds if and only if $D_jv \le w_j$ and $r \le (w_j - D_jv) / \Vert D_j \Vert$. Consequently, for any $v \in \mathcal{V}$, $S(v,r)$ can be reformulated as:
\begin{equation} \notag
	\begin{aligned}
		S(v,r) &= \sum_{j = 1}^{n_c} \epsilon\left(\frac{w_j - D_j v}{\Vert D_j \Vert} - r\right) \\
		&= n_c - \sum_{j = 1}^{n_c} \delta\left(r - \frac{w_j - D_j v}{\Vert D_j \Vert}\right).
	\end{aligned}
\end{equation}
Substituting this into the constraints, the inequality \eqref{4_11_2} becomes:
\begin{equation} \notag
	\max_{v \in \mathcal{V}} \sum_{j = 1}^{n_c} \delta \left(r -\frac{w_j - D_j v}{\Vert D_j \Vert} \right) \le i.
\end{equation}
Thus, the optimization problem \eqref{4_11} is equivalent to:
\begin{subequations} \label{D_4_1}
	\begin{align}
		\max_{r \ge 0} \quad & r, \\
		\text{s.t.} \quad & \max_{v \in \mathcal{V}} \sum_{j = 1}^{n_c} \delta\left(r -\frac{w_j - D_j v}{\Vert D_j \Vert}\right) \le i. \label{D_4_1_con}
	\end{align}
\end{subequations}

Let $h(r) := \max_{v \in \mathcal{V}} \sum_{j = 1}^{n_c} \delta(r - (w_j - D_j v) / \Vert D_j \Vert)$. The feasible set of \eqref{D_4_1} is denoted by $\Gamma_i = \{r \ge 0 \mid h(r) \le i\}$. It can be observed that $\delta(x)$ is non-decreasing and left-continuous. As $h(r)$ is the maximum of finite sums of such functions, $h(r)$ is non-decreasing and left-continuous. To prove that the maximum of \eqref{4_11} is attained, we first show that $\Gamma_i$ is closed. Let $\{r_k\} \subseteq \Gamma_i$ be any convergent sequence with limit $\widetilde{r} = \lim_{k \to \infty} r_k$. We analyze the convergence in two cases:
\begin{enumerate}
	\item[(a)] If there exists a subsequence $\{r_{k_m}\}$ such that $r_{k_m} \ge \widetilde{r}$, we utilize the monotonicity of $h$. Since $r_{k_m} \in \Gamma_i$, we have $h(r_{k_m}) \le i$. By monotonicity, $h(\widetilde{r}) \le h(r_{k_m}) \le i$.
	\item[(b)] If $r_k < \widetilde{r}$ for all sufficiently large $k$, we utilize the left-continuity of $h$. Since $r_k \to \widetilde{r}^-$, we have $h(\widetilde{r}) = \lim_{k \to \infty} h(r_k)$. Given that $h(r_k) \le i$ for all $k$, it follows that $h(\widetilde{r}) \le i$.
\end{enumerate}
In both cases, $\widetilde{r} \in \Gamma_i$, which implies that $\Gamma_i$ is a closed set. Furthermore, since $h(r)$ is non-decreasing and takes integer values, $h(r)$ eventually exceeds $i$ for sufficiently large $r$, making $\Gamma_i$ bounded. Since $\Gamma_i$ is compact and the objective function is continuous, the maximum exists by the Weierstrass extreme value theorem.

(c) Proof of (\ref{4_11_3}): The constraint in \eqref{4_11_1} ensures that the feasible set contains $r=0$, so $\sigma_i \ge 0$. Clearly, the feasible region satisfies $\Gamma_{i} \subseteq \Gamma_{i+1}$ because the condition $h(r) \le i$ implies $h(r) \le i+1$. Consequently, the maximum over a larger set is non-decreasing, i.e., $\sigma_i \le \sigma_{i+1}$.

(d) Proof of (\ref{4_11_5}): It follows directly from monotonicity that
\begin{equation} \notag
	\begin{aligned}
		S(v, r) &\ge S(v, \sigma_i) \\
		&\ge \min_{v \in \mathcal{V}} |\{j \in \mathbb{N}_{[1, n_c]} \ | \ \mathcal{B}(v, \sigma_i) \subseteq \mathcal{V}_j \}| \\
		&\ge n_c - i,
	\end{aligned}
\end{equation}
where the first inequality holds because $\mathcal{B}(v, r) \subseteq \mathcal{B}(v, \sigma_i)$ for any $r \le \sigma_i$.

\subsection{Proof of Lemma \ref{lem:MILP}} \label{app:MILP}
First, we give an alternative computation method of $\sigma_i$ as follows:
\begin{lemma} \label{lem:inf}
	For any $i \in \mathbb{N}_{[1, n_c]}$, $\sigma_i = \overline{\sigma}_i$ where
	\begin{subequations} \label{E_0}
		\begin{align}
			\overline{\sigma}_i = & \inf_{v \in \mathcal{V}, r \ge 0} \ r, \label{E_0_1} \\
			& \text{s.t.} \ \left | \left\{ j \in \mathbb{N}_{[1, n_c]} \ | \ \mathcal{B}(v, r) \subseteq \mathcal{V}_j \right \}\right | < n_c-i. \label{E_0_2}
		\end{align}
	\end{subequations}
	 
\end{lemma}

\begin{pf}
	One the one hand, since $\overline{\sigma}_i$ is the infimum of \eqref{E_0}, for any $\epsilon_1 > 0$, there exist $v_1 \in \mathcal{V}, r_1 \ge 0$ such that 
	\begin{equation}
		\left | \left\{ j \in \mathbb{N}_{[1, n_c]} \ | \ \mathcal{B}(v_1, r_1) \subseteq \mathcal{V}_j \right \}\right | < n_c-i \label{app_E_1}
	\end{equation}
	and $r_1 \le \overline{\sigma}_i + \epsilon_1$. 
	Note that $r_1 > \sigma_i$, otherwise 
	\begin{equation} \notag
		\begin{aligned}
			& \min_{v \in \mathcal{V}} \left | \left\{ j \in \mathbb{N}_{[1, n_c]} \ | \ \mathcal{B}(v, r_1) \subseteq \mathcal{V}_j \right \}\right | \\
			& \ge  \min_{v \in \mathcal{V}} \left | \left\{ j \in \mathbb{N}_{[1, n_c]} \ | \ \mathcal{B}(v, \sigma_i) \subseteq \mathcal{V}_j \right \}\right | \ge n_c-i
		\end{aligned}
	\end{equation}
	which leads to a contradiction with \eqref{app_E_1}. Thus, $\sigma_i \le \overline{\sigma}_{i} + \epsilon_1$ for any $\epsilon_1 > 0$, i.e., $\sigma_i \le \overline{\sigma}_{i}$. 
	
	On the other hand, since $\sigma_{i}$ is the maximum of $\eqref{4_11}$, then for any $\epsilon_2 > 0$, there exists $v_2 \in \mathcal{V}$ such that 
	\begin{equation} \notag
		\left | \left\{ j \in \mathbb{N}_{[1, n_c]} \ | \ \mathcal{B}(v_2, \sigma_i + \epsilon_2) \subseteq \mathcal{V}_j \right \}\right | < n_c-i
	\end{equation}
	which implies that $\overline{\sigma}_i \le \sigma_{i} + \epsilon_2$. Thus, $\overline{\sigma}_i \le \sigma_{i}$. To conclude, $\sigma_{i} = \overline{\sigma}_i$.
\end{pf}

We next show how to reformulate the optimization problem \eqref{E_0} into an MILP. We introduce binary variables $\delta_j \in \{0,1\}, \ j \in \mathbb{N}_{[1, n_c]}$ to represent that  $\mathcal{B}(v,r) \subseteq \mathcal{V}_j$, i.e., $\delta_j = 1$ if $\mathcal{B}(v,r) \subseteq \mathcal{V}_j$ and $\delta_j = 0$, otherwise. Since $\mathcal{B}(v,r) \subseteq \mathcal{V}_j$ if and only if (a) $D_jv\le w_j$
and (b) $r\le |w_j - D_jv| / \Vert D_j \Vert$ (c.f. Remark \ref{rem_formula}), the logical implication of $\delta_j$ can translate into the mixed integer inequality 
\begin{equation} \notag
	\frac{w_j - D_j v}{\Vert D_j \Vert} -M  \delta_j < r \le \frac{w_j - D_j v}{\Vert D_j \Vert} + M (1 - \delta_j),
\end{equation} 
where $M$ is a large number \cite{9928332}. To sum up, the optimization problem in  \eqref{E_0} is reformulated as the following MILP problem:
\begin{equation} \label{E_1}
	\hspace{-0.4cm}
	\begin{aligned}
		\overline{\sigma}_i=& \inf_{v \in \mathcal{V}, r\ge 0} \  r, \\
		\text{s.t.}~ 
		& \frac{w_j - D_j v}{\Vert D_j \Vert} -M  \delta_j < r \le \frac{w_j - D_j v}{\Vert D_j \Vert} + M (1 - \delta_j), \\
		& \sum_{j=1}^{n_c} \delta_j< n_c - i, \ \delta_j \in \{0,1\}, \ j \in \mathbb{N}_{[1, n_c]}.
	\end{aligned}
\end{equation}
By Lemma \ref{lem:inf}, computing $\sigma_i$ amounts to an MILP \eqref{E_1}.

\subsection{Proof of Theorem \ref{thm:adaptive}} \label{app:online_mpc}

(\ref{5_4_2}): Assumption \ref{ass:basic} implies $\Vert x_k \Vert \le c \Vert x_0 \Vert \beta^k$. Thus, $\Vert x_k - x_{k-1} \Vert \le \Vert x_k \Vert + \Vert x_{k-1} \Vert \le c\Vert x_0 \Vert \beta^k(1 + \beta^{-1})$. If $k \ge K_i$, then 
\begin{equation}
	\Vert x_k - x_{k-1}\Vert \le c\Vert x_0 \Vert \beta^{K_i}(1 + \beta^{-1}) = \frac{\sigma_{i}}{\sqrt{1 + \kappa^2}}. \notag
\end{equation}
It follows from Theorem \ref{thm:distance} that $|\mathbb{I}(x_k)| \le n_z + i$. 

(\ref{5_4_3}): It follows from \eqref{indexseti} that
\begin{equation}\notag
	\mathbb{I}(x_k) = \mathbb{A}(x_{k-1}) \bigcup (\mathbb{A}^c(x_{k-1}) - \mathbb{D}(x_k)). 
\end{equation}
Thus, $|\mathbb{I}(x_k)| = 0$ if and only if (i) $\mathbb{A}(x_{k-1}) = \emptyset$ and (ii) $\mathbb{D}(x_k) = \mathbb{N}_{[1, n_c]}$ where $\mathbb{D}(x_k)$ has the form as 
\begin{equation} \label{AA_4_1}
	\{ j \in \mathbb{N}_{[1, n_c]} \ | \ \mathcal{B}(z^*(x_{k-1}), \kappa \Vert x_k - x_{k-1}\Vert) \subseteq \mathcal{Z}_j(x_k) \} 
\end{equation}

We first show that $\mathbb{A}(x_{k-1})= \emptyset$.  Note that the QP(0) of MPC can be trivially solved as $z^*(0)=0$ with $\mathbb{A}^*(0) = \emptyset$. Since the origin is an interior of $\mathcal{X}, \mathcal{U} \text{  and } \mathcal{X}_N$, all constraints in the $\text{QP}(0)$ are inactive, i.e., $G  0 < w + S z^*(0)$. Thus, $w > 0$. Let $\widehat{x} = \textbf{0}$ and $\mathbb{I}(x) = \mathbb{A}(x_{k-1})$ in \eqref{ball}. It follows that  $z^*(x_{k-1}) \in \mathcal{B}(\textbf{0}, \kappa \Vert x_{k-1} \Vert)$. Then, 
\begin{equation} \notag
	\begin{aligned} 
		& G_jz^*(x_{k-1}) - S_j x_{k-1} - w_j \\ 
		& \le \kappa \Vert G_j \Vert \Vert x_{k-1} \Vert + \Vert S_j \Vert \Vert x_{k-1} \Vert - w_j \\
		& \le c\Vert x_0 \Vert(\kappa \Vert G_j \Vert + \Vert S_j \Vert) \beta^{k-1} - w_j \\
		& \le c\Vert x_0 \Vert(\kappa \Vert G_j \Vert + \Vert S_j \Vert) \beta^{\widehat{K}_1 - 1} - w_j.
	\end{aligned}
\end{equation}
where the second inequality follows from Assumption \ref{ass:basic} and the third from $k \ge \widehat{K} \ge \widehat{K}_1$. Substituting the definition of $\widehat{K}_1$,  $G_jz^*(x_{k-1}) - S_j x_{k-1} - w_j < 0$ for any $j \in \mathbb{N}_{[1, n_c]}$, i.e., $\mathbb{A}(x_{k-1}) = \emptyset$. 

We next show that $\mathbb{D}(x_k) = \mathbb{N}_{[1, n_c]}$. According to \eqref{AA_4_1}, $j \in \mathbb{D}(x_k)$ if and only if (i) $G_j z^*(x_{k-1}) \le S_j x_k + w_k$ and (ii) 
$|w_j + S_j x_k - G_j z^*(x_{k-1})|/\Vert G_j \Vert  \ge \kappa \Vert x_k - x_{k-1} \Vert$ (c.f. Remark \ref{rem_formula}). The necessary and sufficient condition can be merged as 
\begin{equation} \label{AA_4_2}
	\kappa \Vert G_j \Vert \Vert x_k - x_{k-1} \Vert \le w_j + S_j x_k - G_j z^*(x_{k-1})
\end{equation}
Since $\mathbb{A}(x_{k-1})= \emptyset$, $z^*(x_{k-1}) = -H^{-1}F^Tx_{k-1}$. Then, $G_j z^*(x_{k-1}) \le \Vert  G_j H^{-1}F^T\Vert \Vert x_{k-1}\Vert $. In fact, $S_j x_k \ge -\Vert S_j \Vert \Vert x_k \Vert$ and $\Vert x_k - x_{k-1} \Vert \le \Vert x_k \Vert + \Vert x_{k-1} \Vert $. Hence, the sufficient condition of \eqref{AA_4_2} can be obtained as follows:
\begin{equation} \notag
	\begin{aligned}
		 &  w_j \ge \kappa \Vert G_j \Vert (\Vert x_k \Vert + \Vert x_{k-1} \Vert) + \Vert S_j \Vert \Vert x_k \Vert + \Vert  G_j H^{-1}F^T\Vert  \\
		 & \times \Vert x_{k-1}\Vert \ge c\Vert x_0 \Vert \beta^k ( \kappa \Vert G_j \Vert (1 + \beta^{-1}) + \Vert S_j \Vert  \\
		& + \Vert  G_j H^{-1}F^T\Vert \beta^{-1}) \ge 0.
	\end{aligned}
\end{equation}
where the second inequality follows from Assumption \ref{ass:basic}. Since $w> 0 $, the sufficient condition holds for any $j \in \mathbb{N}_{[1, n_c]}$. As a result, $\mathbb{D}(x_k) = \mathbb{N}_{[1, n_c]}$. To sum up, $| \mathbb{I}(x_k)| = 0$. 	

\bibliographystyle{IEEEtran}      
\bibliography{tac}           

 \begin{IEEEbiography}[{\includegraphics[width=1in,height=1.25in,clip,keepaspectratio]{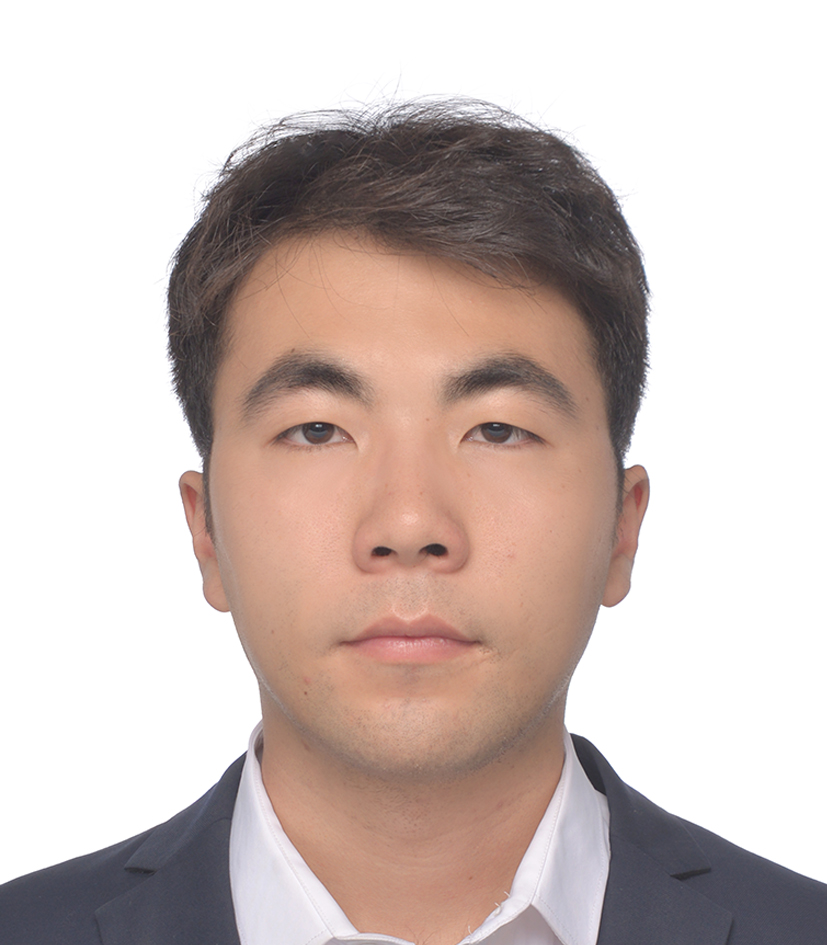}}]{Zhinan Hou} received the B.S. degree in Control Science and Engineering at the Department of Automation, Tsinghua University, Beijing, China, in 2022, where he is currently working toward the Ph.D. degree. His research interests include multiparametric quadratic program, integer programming, reinforcement learning.
 \end{IEEEbiography}
 
 \vspace{-405pt}
 
 \begin{IEEEbiography}
 	[{\includegraphics[width=1in,height=1.25in,clip,keepaspectratio]{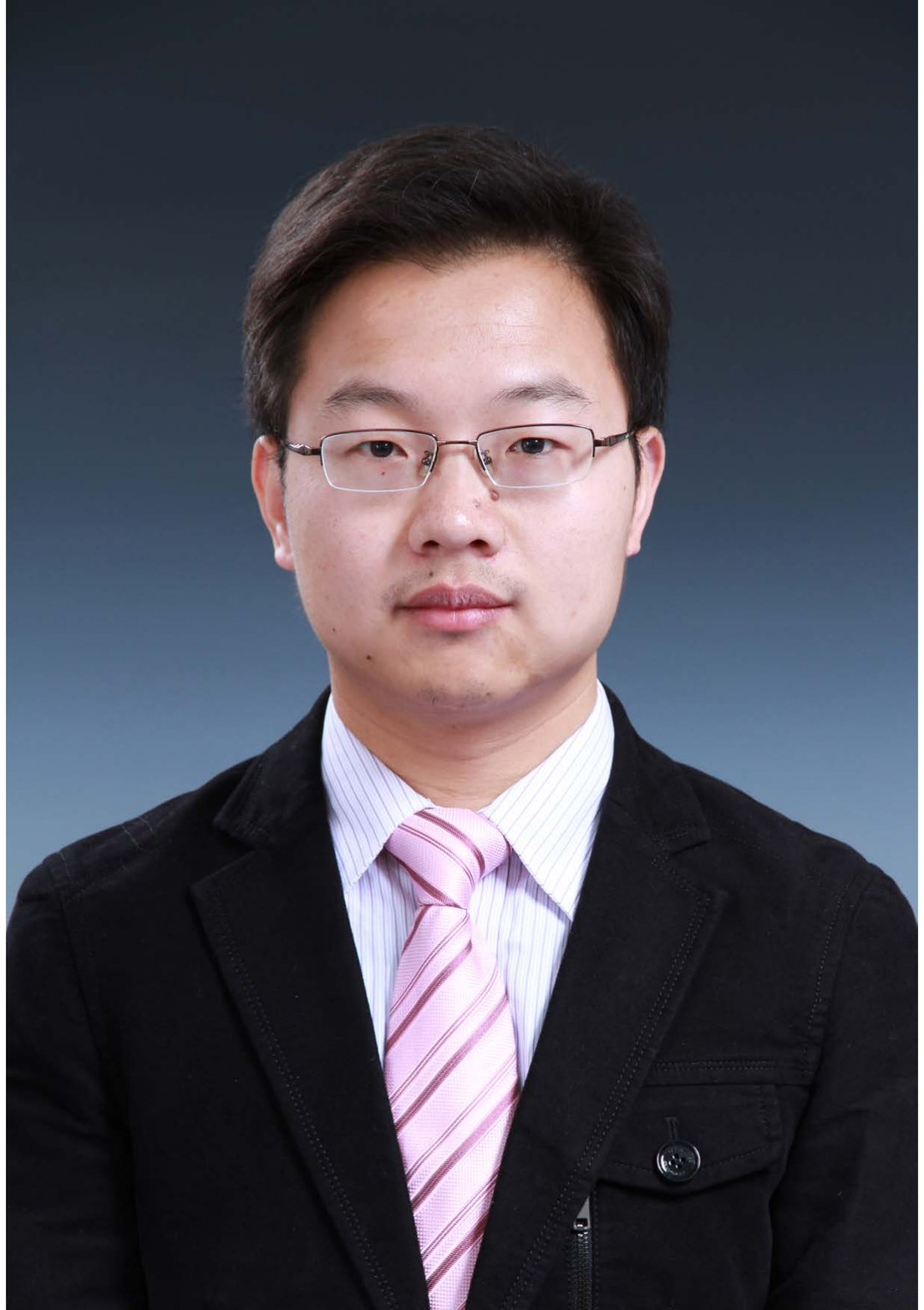}}]
 	{Keyou You} received the B.S. degree in Statistical Science from Sun Yat-sen University, Guangzhou, China, in 2007 and the Ph.D. degree in Electrical and Electronic Engineering from Nanyang Technological University (NTU), Singapore, in 2012. After briefly working as a Research Fellow at NTU, he joined Tsinghua University in Beijing, China where he is now a Full Professor in the Department of Automation. He held visiting positions at Politecnico di Torino, Hong Kong University of Science and Technology, University of Melbourne and etc.
 	
 	Prof. You's research interests focus on the intersections between control, optimization and learning as well as their applications in autonomous systems. He received the Guan Zhaozhi award at the 29th Chinese Control Conference in 2010 and the ACA (Asian Control Association) Temasek Young Educator Award in 2019. He received the National Science Funds for Excellent Young Scholars in 2017, and for Distinguished Young Scholars in 2023. Currently, he is an Associate Editor for {\em Automatica}, {\em IEEE Transactions on Automatic Control} and  {\em IEEE Transactions on Control of Network Systems}.
 \end{IEEEbiography}

	
\end{document}